\documentclass[11pt]{amsart}
\usepackage {amsmath}
\usepackage {amsthm}
\usepackage {enumerate}
\usepackage{amssymb}
\usepackage{amsfonts}
\usepackage{anysize}
\usepackage{setspace}
\usepackage[all]{xy}

\marginsize{3.15cm}{3.15cm}{1.6cm}{3cm}

\newcommand{\N}{\mathbb N}

\newcommand{\R}{\mathbb R}

\numberwithin{equation}{section}

\theoremstyle{plain}

\newtheorem{theorem}[equation]{Theorem}
\newtheorem{proposition}[equation]{Proposition}
\newtheorem{corollary}[equation]{Corollary}
\newtheorem{lemma}[equation]{Lemma}
\newtheorem{question}[equation]{Question}
\newtheorem{problem}[equation]{Problem}

\theoremstyle{definition}
\newtheorem{definition}[equation]{Definition}

\theoremstyle{definition}

\theoremstyle{definition}

\theoremstyle{definition}
\newtheorem{example}[equation]{Example}

\theoremstyle{definition}

\author{{\bfseries S. Garc\'ia-Ferreira}}
\address{Centro de Ciencias Matem\'aticas\\
         Universidad Nacional Aut\'onoma de M\'exico\\
				 Campus Morelia\\
         Apartado Postal 61-3, Santa Mar\'ia, 58089, Morelia, Michoac\'an, M\'exico.}
\email{sgarcia@matmor.unam.mx}

\title{\scshape\bfseries families of continuous retractions and function spaces}

\author{{\bfseries R. Rojas-Hern\'andez}}
\email{satzchen@yahoo.com.mx}



\date{}

\thanks{Research of the first-named author was supported
by  CONACYT grant no. 176202 and PAPIIT grant no. IN-101911.}

\begin{document}

\begin{abstract} In this article, we mainly study certain families of continuous retractions ($r$-skeletons) having certain rich properties. By using  monotonically retractable spaces we solve a question posed by R. Z. Buzyakova in  \cite{buz} concerning the Alexandroff duplicate of a space. Certainly, it is shown that if the space  $X$ has a  full $r$-skeleton, then its Alexandroff duplicate  also has a full  $r$-skeleton and, in a very similar way, it is proved that the Alexandroff duplicate of a monotonically retractable space is monotonically retractable. The notion of $q$-skeleton is introduced and it is shown that every compact subspace of $C_p(X)$ is Corson when $X$ has a full $q$-skeleton. The notion of strong $r$-skeleton is also introduced to answer  a question suggested by F. Casarrubias-Segura and R. Rojas-Hern\'andez in their paper \cite{cas-rjs} by  establishing  that a space $X$ is monotonically Sokolov iff it is monotonically $\omega$-monolithic and has a strong $r$-skeleton. The techniques used here allow us to give a  topological  proof of a result of  I. Bandlow \cite{ban} who used elementary submodels and uniform spaces.
\end{abstract}

\maketitle

\section{Preliminaries and Introduction}

Our spaces will be Tychonoff (completely regular and Hausdorff). The set of natural numbers will be denoted by $\mathbb{N}$ and $\omega$ will stand for the first infinite cardinal number. Given an infinite set  $X$, $\mathcal{P}(X)$ is the power set of $X$, the symbol $[X]^{\leq \omega}$ will denote the set of all countable  subsets of $X$ and the meaning of $[X]^{< \omega}$ should be clear.  The  real line with the usual order topology will be denoted by  $\R$ and $\mathcal{B}(\R)$ will stand for a countable fixed base for the topology of  $\R$. For a space $X$, $C_p(X)$ will be the set $C(X)$ of all real-valued continuous functions on $X$ equipped with the topology of pointwise convergence and $exp(X)$ will be the family of all non-empty closed subsets of $X$. For a continuous map $f : X \to Y$ we denote by $f^{\ast} : C_p(Y) \to C_p(X)$ the dual map of $f$ given by $f^{\ast}(g) = g \circ f$ for all $g \in C_p(Y)$. If $Y \subseteq X$, then we denote by $\pi_Y : C_p(X) \to C_p(Y)$ the function which restricts any map in $C_p(X)$ to $Y$. For a nonempty set $A \subseteq C_p(X)$ the map $\Delta_{A} : X \to \R^A$ is called the {\it diagonal map} of $A$. A surjective map $f : X \to Y$ is called {\it ${\mathbb R}$-quotient} if, for every function  $g : Y \to {\mathbb R}$ the continuity of the composition $g \circ f$ implies the continuity of $g$. A continuous surjection  will be referred as a  {\it condensation}. A space $X$ is called {\it cosmic} if it has a countable network. For  $\mathcal{N} \subseteq \mathcal{P}(X)$ and $f : X \to Y$, we say that $\mathcal{N}$ is a {\it network of $f$} if for every $x \in X$ and each open set $U$ in $Y$ with $f(x) \in U$ there is $N \in \mathcal{N}$ such that $x \in  N$ and $f (N) \subseteq U$. Given a set $A \subseteq X$,  a family $\mathcal{N} \subseteq exp(X)$ is said to be an {\it external network} of $A$ in $X$ if for each $x \in A$ and each open set $U$ with $x \in U$ there exists $N \in \mathcal{N}$ such that $x \in N \subseteq U$. All topological notions whose definitions are not stated explicitly here should be understood as  in  \cite{ark} and \cite{ltka} .

\medskip

The families of  continuous retractions often appear in Functional Analysis and in General To\-po\-lo\-gy: For instance,  the   projectional resolutions of the identity (see \cite{fa}) and the inverse limits of retractions (see \cite{kbs-mich}).  One more example of this kind of families is the following notion of $r$-skeleton:

\begin{definition}{\bf\cite{kbs-mich}}\label{DRS}
An {\it $r$-skeleton} in a space $X$ is a family of continuous retractions $\{r_s : s \in \Gamma\}$, indexed by an up-directed and $\sigma$-complete partially ordered set $\Gamma$, such that:
\begin{enumerate}[(i)]
	\item  $r_s(X)$ is cosmic for each $s \in \Gamma$,
  \item  $r_s = r_s \circ r_t = r_t \circ r_s$ whenever $s \leq t$,
  \item if $\{s_n : n \in \mathbb{N}\} \subseteq \Gamma$, $s_n \leq s_{n+1}$ for each $n \in \mathbb{N}$ and $t = \sup_{n \in \mathbb{N}} s_n$; then $r_t(x) =  \lim_{n \to \infty} r_{s_n}(x)$ for each $x \in X$, and

  \item $x = \lim_{s \in \Gamma} r_s(x)$ for every $x \in X$.
\end{enumerate}
If $X = \bigcup_{s \in \Gamma} r_s(X)$, then we say that $\{r_s : s \in \Gamma\}$ is a {\it full $r$-skeleton}.
\end{definition}

W. Kubi\'s and H. Michalewski \cite{kbs-mich} used the $r$-skeletons to characterize the Valdivia compact spaces: A compact space is Valdivia  iff it has a commutative $r$-skeleton. In the paper \cite{mrk2}, the author applied the  $r$-skeletons to characterize Corson compact spaces: A compact space is Corson if and only if it has a full $r$-skeleton. The $r$-skeletons have also been a very important tool in the study of certain topological properties. For instance,
M. C\'uth and O. Kalenda \cite{mrk-kld} proved that a countably compact space is monotonically retractable iff it has a full $r$-skeleton. In this way, they answered a question posed in \cite{rjs-tka}, by showing that a compact space is monotonically retractable if and only if it is Corson.
On the other hand,  I. Bandlow  \cite{ban2} gave a characterization of the  Corson-compact spaces by means of countable elementary substructures. This characterization certainly motivates  the consideration of  the $r$-skeletons. Three years later, by using also elementary submodels  and uniform structures,  Bandlow \cite{ban}  solved  a problem posed by  A. V. Arhangelskii \cite[Prob. IV.3.16]{ark} by establishing  that a compact space $X$ is Corson when $C_p(X)$ is  a continuous image of a closed subspace of the product $L_\kappa^\omega \times K$, where $K$ is an arbitrary compact space and $L_\kappa$ denotes the Lindel\"of extension by one point of the discrete space of infinite cardinality $\kappa$.

\medskip

The second section is devoted to show basic properties of $\omega$-monotone maps that will be used in the subsequent sections.
Our main result in the third section  is to  prove that if a space either has a full $r$-skeleton or it is monotonically retractable, then its Alexandroff duplicate has the same property. As a consequence, the Alexandroff duplicate of a Corson compact is also Corson compact. Besides, we show that  the Alexandroff duplicate of a $\Sigma$-product of real lines is monotonically retractable, this solves positively \cite[Question 3.16]{buz} (see Corolary \ref{CSPB}).
In the fourth section, we introduced the notion of $q$-skeleton in parallel sense to the notion of $r$-skeleton. We show that every compact subspace of $C_p(X)$ has a full $r$-skeleton (i.e., it is Corson) whenever $X$ has a full $q$-skeleton. We also prove that Lindel\"of $\Sigma$-spaces and pseudocompact spaces have a full $q$-skeleton.
Monotonically Sokolov spaces where studied in \cite{rjs-tka}, where a two-way $C_p$-duality between monotonically retractable and monotonically Sokolov spaces was established. In the fifth section, we introduce the notion of strong $r$-skeleton and show that a space $X$ is monotonically Sokolov iff it is monotonically $\omega$-monolithic and has a strong $r$-skeleton. This result answers a question in \cite{cas-rjs} and clarifies the relation between monotonically retractable spaces, monotonically Sokolov spaces, $q$-skeletons and $r$-skeletons.
Finally, we apply $q$-skeletons, in the section fifth, to provide a  topological proof of Bandlow's result quoted above.

\section{Monotone assignments and monotone properties}

We start this section with the description of notion a $\omega$-monotone map.

\begin{definition}
Let  $\Gamma$  be an up-directed and $\sigma$-complete partially ordered set and $Y$ a set. A function  $\phi : \Gamma \to [Y]^{\leq\omega}$  is called
  $\omega$-{\it monotone} provided that:
\begin{enumerate}[(a)]
	\item if $s,t \in \Gamma$ and $s \leq t$, then $\phi(s) \subseteq \phi(t)$; and
  \item if $\{s_n : n \in \mathbb{N}\} \subseteq \Gamma$, $s_n \leq s_{n+1}$ for each $n \in \mathbb{N}$, and $t = \sup_{n \in \mathbb{N}} s_n$, then
  $$\phi(t) =  \bigcup_{n \in \mathbb{N}} \phi(s_n).$$
\end{enumerate}
\end{definition}

We are mainly interested in the following topological properties that somehow involve monotone maps.

\medskip

For the case when $\Gamma = [X]^{\leq\omega}$, where $X$ is a nonempty set, the order involved will be always  the containment.

\begin{definition}
 Let $X$ be a space. We say that $X$ is:
\begin{itemize}
\item \cite{atw} \textit{monotonically $\omega$-monolithic}  if there exits $\mathcal{N} : [X]^{\leq\omega} \to [\mathcal{P}(X)]^{\leq\omega}$ such that $\mathcal{N}(A)$ is an external network of $\textnormal{cl}_X(A)$ for each $A \in [X]^{\leq\omega}$ and $\mathcal{N}$ is $\omega$-monotone.

\item \cite{rjs2} \textit{monotonically $\omega$-stable} if there exits an $\omega$-monotone map $\mathcal{N} : [C_p(X)]^{\leq\omega} \to [\mathcal{P}(X)]^{\leq\omega}$ such that $\mathcal{N}(A)$ is a network of $\Delta_{\textnormal{cl}(A)}$, for all $A \in [C_p(X)]^{\leq\omega}$.

\item \cite{rjs1} \textit{monotonically retractable} if we can assign to each $A \in [X]^{\leq\omega}$ a continuous retraction
$r_A : X \to X$ and a family $\mathcal{N}(A) \in [\mathcal{P}(X)]^{\leq\omega}$ such that $A \subseteq r_A(X)$, $\mathcal{N}(A)$ is a network of $r_A$ and $\mathcal{N}$ is $\omega$-monotone.

\item \cite{rjs-tka} \textit{monotonically Sokolov} if we can assign to every $\mathcal{F} \in [exp(X)]^{\leq\omega}$ a continuous retraction $r_\mathcal{F} : X \to X$ and a family $\mathcal{N}(\mathcal{F}) \in [\mathcal{P}(X)]^{\leq\omega}$ such that $r_\mathcal{F}(F) \subseteq F$ for each $F \in \mathcal{F}$, $\mathcal{N}(\mathcal{F})$ is an external network of $r_\mathcal{F}(X)$ and $\mathcal{N}$ is $\omega$-monotone.
\end{itemize}
\end{definition}

The next $\omega$-monotone map will be essential in some proofs.

\begin{example}\label{EOW}
 Given $N_1,\ldots,N_n \in \mathcal{P}(X)$ and $U_1,\ldots,U_n \in \mathcal{P}(\R)$, put
 $$[N_1,\ldots, N_n;U_1,\ldots, U_n] = \{f \in C_p(X) : \forall i \leq n \big(f(N_i) \subseteq U_i\big) \}.$$
 Consider the map $\mathcal{W} : [\mathcal{P}(X)]^{\leq\omega} \to [\mathcal{P}(C_p(X))]^{\leq\omega}$ defined by
 $$\mathcal{W}(\mathcal{N}) = \{[N_1,\ldots, N_n; B_1,\ldots, B_n] :  n \in \mathbb{N} \ \text{and} \ \forall i \leq n  \big(N_i \in \mathcal{N} \wedge B_i \in \mathcal{B}(\R) \big) \},$$
 for each $\mathcal{N} \in [\mathcal{P}(X)]^{\leq\omega}$. It is straightforward to verify that $\mathcal{W}$ is $\omega$-monotone. For each $D \in [X]^{\leq\omega}$ let $\mathcal{W}_0(D) = \mathcal{W}(\{\{x\} : x \in D\})$ which is the family of all canonical open sets in $C_p(X)$ with support in $D$. It is easy to see that the map $\mathcal{W}_0 : [X]^{\leq\omega} \to [\mathcal{P}(C_p(X))]^{\leq\omega}$ is also $\omega$-monotone.
\end{example}

In what follows we shall frequently and without reference the following two easy facts.

\begin{itemize}
\item If $\phi : [X]^{\leq\omega} \to [Y]^{\leq\omega}$ and $\psi : [Y]^{\leq\omega} \to [Z]^{\leq\omega}$ are $\omega$-monotone, then $\psi \circ \phi$ is $\omega$-monotone.

\item If $\phi : [X]^{\leq\omega} \to [Y]^{\leq\omega}$, then $\phi$ is $\omega$-monotone if and only if there exists a map $\varphi : [X]^{<\omega} \to [Y]^{\leq\omega}$ such that $\phi(A) = \bigcup\{\varphi(F) : F \in [A]^{<\omega}\}$.
\end{itemize}

\medskip

Now, we shall prove two  basic results.

\begin{definition}\label{DCOM}
Let $\phi : [X]^{\leq\omega} \to [X]^{\leq\omega}$ be a function. For a set  $A \in [X]^{\leq\omega}$, the  {\it closure of $A$ under $\phi$} is the set
$$
\overline{\phi}(A) = \bigcap\left\{B \in [X]^{\leq\omega} : A \subseteq B \ \text{and} \ \phi(B) \subseteq B\right\}.
$$
\end{definition}

The next lemma give us a method to generate $\omega$-monotone assignments with some nice properties which are very useful for our purposes.

\begin{lemma}\label{LCOM}
 Let $\phi : [X]^{\leq\omega} \to [X]^{\leq\omega}$ be an $\omega$-monotone map. Then for each $A \in [X]^{\leq\omega}$ the set $\overline{\phi}(A)$ satisfies that  $A \subseteq \overline{\phi}(A) \in [X]^{\leq\omega}$ and the function $\overline{\phi}$  is $\omega$-monotone.
\end{lemma}

\proof
We are going to construct $\overline{\phi}(A)$ by a recursive process. Let $\phi_0(A) = A$ and $\phi_{n+1}(A) = \phi_n(A) \cup \phi(\phi_n(A))$ for each $n \in \mathbb{N}$. Using induction it is easy to verify that $\phi_n(A) \in [X]^{\leq\omega}$ and $\phi_n(A) \subseteq \overline{\phi}(A)$ for each $n \in \mathbb{N}$. Hence, $\bigcup_{n \in \mathbb{N}}\phi_n(A) \subseteq \overline{\phi}(A)$. On the other hand $A \subseteq \bigcup_{n \in \mathbb{N}}\phi_n(A)$ and since $\phi$ is $\omega$-monotone, we obtain that $\phi(\bigcup_{n \in \mathbb{N}}\phi_n(A))$ $\subseteq \bigcup_{n \in \mathbb{N}}\phi_n(A)$. It follows from the definition that $\overline{\phi}(A) \subseteq \bigcup_{n \in \mathbb{N}}\phi_n(A)$. So, $\overline{\phi}(A) = \bigcup_{n \in \mathbb{N}}\phi_n(A)$ and $A \subseteq \overline{\phi}(A) \in [X]^{\leq\omega}$. By applying induction again, we can verify that $\phi_n$ is $\omega$-monotone for each $n \in \mathbb{N}$.  As a consequence, we have the following:
\begin{enumerate}[(a)]
	\item If $A,B \in [X]^{\leq\omega}$ and $A \subseteq B$, then $\overline{\phi}(A) = \bigcup_{n \in \mathbb{N}}\phi_n(A) \subseteq \bigcup_{n \in \mathbb{N}}\phi_n(B) = \overline{\phi}(B).$
  \item If $\{A_m : m \in \mathbb{N}\} \subseteq [X]^{\leq\omega}$, $A_m \leq A_{m + 1}$ for each $m \in \mathbb{N}$ and $A = \bigcup_{m \in \mathbb{N}}A_m$; then $\overline{\phi}(A) = \bigcup_{n \in \mathbb{N}}\phi_n(A) = \bigcup_{n \in \mathbb{N}}\bigcup_{m \in \mathbb{N}}\phi_n(A_m) = \bigcup_{m \in \mathbb{N}}\bigcup_{n \in \mathbb{N}}\phi_n(A_m) = \bigcup_{m \in \mathbb{N}} \overline{\phi}(A_m).$
\end{enumerate}
Therefore, $\overline{\phi}$ is $\omega$-monotone.
\endproof

\begin{lemma}\label{LAMG}
Let $X$ be a set and let $\Gamma$ be an up-directed and $\sigma$-complete partially ordered set. Suppose that for each $x \in X$ we have assigned $s_x \in \Gamma$. Then there exists a function  $\gamma : [X]^{\leq\omega} \to \Gamma$ such that
\begin{enumerate}[(i)]
  \item $\gamma(\{x\}) \geq s_x$ for each $x \in X$;
	\item if $A \subseteq B$, then $\gamma(A) \leq \gamma(B)$; and
	\item if $\{A_n : n \in \mathbb{N}\} \subseteq [X]^{\leq\omega}$, $A_n \leq A_{n+1}$ for each $n \in \mathbb{N}$ and $A = \bigcup_{n \in \mathbb{N}}A_n$; then  $\gamma(A) =  \sup_{n \in \mathbb{N}} \gamma(A_n)$.
\end{enumerate}
\end{lemma}

\proof
First we define $\gamma(F)$ for $F \in [X]^{<\omega}$. The point $\gamma(F)$ will be defined by induction on the cardinality of $F$ as follows: $\gamma(\emptyset) = s_0$ for some arbitrary $s_0 \in \Gamma$, and for a nonempty $F \in [X]^{<\omega} $ we let $\gamma(F)$ be an upper bound of $\{\gamma(G) : G \subseteq F, G \not= F\} \cup \{s_x : x \in F\}$.  For the general case, we define $\gamma(A) = \sup\{\gamma(F) : F \in [A]^{<\omega}\}$ for every $A \in [X]^{\leq\omega}$. First note that $\gamma(A)$ is well defined and it is not hard to see that $\gamma$ satisfies (i) and (ii). To verify (iii) we assume that $A = \bigcup_{n \in \mathbb{N}} A_n \in [X]^{\leq\omega}$ where $A_n \subseteq A_{n+1}$ for each $n \in \mathbb{N}$. Because of (ii) we have $\gamma(A_n) \leq \gamma(A_{n+1}) \leq \gamma(A)$ for each $n \in \mathbb{N}$, and so $\sup_{n \in \mathbb{N}} \gamma(A_n) \leq \gamma(A)$. On the other hand, if $F \in [A]^{<\omega}$, then there exists $m \in \mathbb{N}$ such that $F \in [A_m]^{<\omega}$ and hence $\gamma(F) \leq \gamma (A_m) \leq \sup_{n \in \mathbb{N}}\gamma(A_n)$. It then follows that $\gamma(A) \leq \sup_{n \in \mathbb{N}}\gamma(A_n)$. Therefore, $\gamma(A) = \sup_{n \in \mathbb{N}}\gamma(A_n)$.
\endproof

\begin{proposition}\label{PEM}
If the space $X$ has a full $r$-skeleton, then $X$ has a full $r$-skeleton $\{r_A : A \in [X]^{\leq\omega}\}$ such that for each $A \in [X]^{\leq\omega}$ we have that $r_A(x) = x$  for every $x \in A$.
\end{proposition}

\proof
Assume that $\{r_s : s \in \Gamma\}$ is a full $r$-skeleton in $X$. Given $x \in X$ fix $s_x \in \Gamma$ such that $r_{s_x}(x) = x$. Consider $\gamma : [X]^{\leq\omega} \to \Gamma$ as in Lemma \ref{LAMG}. For $A \in [X]^{\leq\omega}$ we define $r_A = r_{\gamma(A)}$. We left the reader to check that $\{r_A : A \in [X]^{\leq\omega}\}$ is a full $r$-skeleton which satisfies the desired conditions.
\endproof

\section{The Alexandroff duplicate of a space}

Let us remind the definition of the \textit{Alexandroff duplicate} $AD(X)$ of a space $X$: It is  the space $X \times \{0,1\}$ with the topology in which all points of $X \times \{1\}$ are isolated, and basic neighborhoods of points $(x,0)$ are of the form $(U \times \{0,1\}) \setminus \{(x,1)\}$ where $U$ is a neighborhood of $x$ in $X$. We denote by $\pi$ the projection from $AD(X)$ onto $X$.
Our task in this section is to  prove that if a space either has a full $r$-skeleton or it is monotonically retractable, then its Alexandroff duplicate has the same property.

\begin{theorem}\label{ADFRS}
If $X$ has a full $r$-skeleton, then $AD(X)$ also has a full $r$-skeleton.
\end{theorem}

\proof
Consider a full $r$-skeleton $\{r_A : A \in [X]^{\leq\omega}\}$ on $X$ as in Proposition \ref{PEM}. For each $A \in [AD(X)]^{\leq\omega}$ we define the map $\hat{r}_A : AD(X) \to AD(X)$ as follows:
\[
\hat{r}_A((x,i)) =
\begin{cases}
(r_{\pi(A)}(x),i) & \textrm{ if $x \in \pi(A)$}\\
(r_{\pi(A)}(x),0) & \textrm{ if $x \in X \setminus \pi(A)$}, \\
\end{cases}
\]
for $x \in X$ and $i \in \{0,1\}$. For each $A \in [AD(X)]^{\leq\omega}$ observe that $\hat{r}_A$ is a continuous retraction and $\hat{r}_A((x,i)) = (r_{\pi(A)}(x),i) = (x,i)$ for each $x \in \pi(A)$ and $i \in \{0,1\}$. Let  $\Gamma := [AD(X)]^{\leq\omega}$. We shall prove that $\{\hat{r}_A : A \in \Gamma\}$ is a full $r$-skeleton on $AD(X)$.\medskip

(i) Given $A \in \Gamma$, since $A$ is countable and $r_{\pi(A)}(X)$ is cosmic, the space $\hat{r}_A(AD(X)) = \big(r_{\pi(A)}(X) \times \{0\}\big) \cup \big(\pi(A) \times\{1\}\big)$ is also cosmic.\medskip

(ii) Assume that $A,B \in \Gamma$ and $A \subseteq B$. Note that $\pi(A) \subseteq \pi(B)$. Choose $(x,i) \in AD(X)$. We shall verify that $\hat{r}_A((x,i)) = \hat{r}_A \circ \hat{r}_B ((x,i)) = \hat{r}_B \circ \hat{r}_A ((x,i))$. Consider the following tree cases.\medskip

\textit{Case 1.}   $x \in X \setminus \pi(B) \subseteq X \setminus \pi(A)$. Then we have
$\hat{r}_A(\hat{r}_B((x,i))) = \hat{r}_A((r_{\pi(B)}(x),0)) = (r_{\pi(A)}(r_{\pi(B)}(x)),0)  = (r_{\pi(A)}(x),0) = \hat{r}_A((x,i))$
and
$\hat{r}_B(\hat{r}_A((x,i))) = \hat{r}_B((r_{\pi(A)}(x),0)) = (r_{\pi(B)}(r_{\pi(A)}(x)),0)  = (r_{\pi(A)}(x),0) = \hat{r}_A((x,i))$.\medskip

\textit{Case 2.}  $x \in \pi(B) \setminus  \pi(A)$. In this case, we obtain that  $\hat{r}_A(\hat{r}_B((x,i))) =  \hat{r}_A((x,i))$ and
$\hat{r}_B(\hat{r}_A((x,i))) =  \hat{r}_B((r_{\pi(A)}(x),0))  = (r_{\pi(B)}(r_{\pi(A)}(x)),0) =  (r_{\pi(A)}(x),0)  =  \hat{r}_A((x,i)).$\medskip

\textit{Case 3.}  $x \in \pi(A) \subseteq \pi(B)$. So,
$$
\hat{r}_A(\hat{r}_B((x,i))) =  \hat{r}_A((x,i))  = (x,i) = \hat{r}_B((x,i)) = \hat{r}_B(\hat{r}_A((x,i))).
$$

(iii) Choose $\{A_n : n \in \mathbb{N}\} \subseteq \Gamma$ with $A_n \subseteq A_{n+1}$ for each $n \in \mathbb{N}$ and let $A = \bigcup_{n \in \mathbb{N}} A_n$. We will prove that $\hat{r}_A((x,i)) = \lim_{n \to \infty} r_{A_n}((x,i))$ for every $(x,i) \in AD(X)$. Let $(x,i) \in AD(X)$. Note that $\pi(A) = \bigcup_{n \in \mathbb{N}}\pi(A_n)$ and hence $r_{\pi(A)}(x) = \lim_{n \to \infty} r_{\pi(A_n)}(x)$. If $x \in X \setminus \pi(A)$, then $x \in X \setminus \pi(A_n)$ for any $n \in \mathbb{N}$, and so $\hat{r}_A((x,i)) = (r_{\pi(A)}(x),0) = (\lim_{n \to \infty} r_{\pi(A_n)}(x),0) = \lim_{n \to \infty} (r_{\pi(A_n)}(x),0) = \lim_{n \to \infty} \hat{r}_{A_n}((x,i))$. If $x \in \pi(A)$, then $x \in \pi(A_n)$  when $n > m$ for some $m \in \N$. Therefore, $\lim_{n \to \infty} r_{A_n}((x,i)) = \lim_{n \to \infty}(x,i) = (x,i) = r_{\pi(A)}((x,i))$.\medskip

(iv) Given $(x,i) \in AD(X)$ there exists $\{(x,i)\} \in \Gamma$ such that if $A \in \Gamma$ and $\{(x,i)\} \subseteq A$, then $x \in \pi(A)$ and so $\hat{r}_A((x,i)) = (x,i)$. It follows that $(x,i) = \lim_{A \in \Gamma}\hat{r}_A((x,i))$.\\

Finally, it is clear that $AD(X) = \bigcup_{A \in \Gamma}\hat{r}_A(AD(X))$.
\endproof

\begin{corollary}\cite{kld}
If $X$ is Corson compact, then $AD(X)$ is also Corson compact.
\end{corollary}

\proof
Let $X$ be a Corson compact space. Because of the characterization of Corson compact spaces obtained in \cite[Theorem 3.11]{mrk2}, the space $X$ has a full $r$-skeleton. Then we can apply Theorem \ref{ADFRS} to see that $AD(X)$ also has a full $r$-skeleton. As $AD(X)$ is a compact space,  applying  \cite[Theorem 3.11]{mrk2} again, we obtain that $AD(X)$ is a Corson compact space.
\endproof

In particular, we have the following consequence.

\begin{corollary}\label{CPB}
If $X$ is Corson compact, then $C_p(AD(X))$ is Lindel\"of.
\end{corollary}

In \cite{buz} Buzyakova asked whether or not the function space over the Alexandroff duplicate of a $\Sigma$-product of real lines has the Lindel\"of property. In order to get a positive answer to this question (in Corollary \ref{CSPB}), we shall prove that monotone retractability is preserved under Alexandroff duplicates.

\begin{theorem}\label{ADMR}
If $X$ is monotonically retractable, then $AD(X)$ is also monotonically retractable.
\end{theorem}

\proof
Let $X$ be a monotonically retractable space. Hence, by definition, we  assign to each $A \in [X]^{\leq\omega}$ a continuous retraction $r_A : X \to X$ and a family $\mathcal{N}(A) \in [\mathcal{P}(X)]^{\leq\omega}$ such that $A \subseteq r_A(X)$, $\mathcal{N}(A)$ is a network of $r_A$ and $\mathcal{N}$ is $\omega$-monotone. For each $A \in [AD(X)]^{\leq\omega}$ consider the continuous retract $\hat{r}_A : AD(X) \to AD(X)$ as in the proof of Theorem \ref{ADFRS}. As above, for each $A \in [AD(X)]^{\leq\omega}$ we know that $\hat{r}_A((x,i)) = (r_{\pi(A)}(x),i) = (x,i)$ for each $x \in \pi(A)$ and $i \in \{0,1\}$. Now, for each $A \in [AD(X)]^{\leq\omega}$ we define the set
$$
\hat{\mathcal{N}}(A) = \{(N \setminus F) \times \{i\} : N \in \mathcal{N}(\pi(A)), F \in[\pi(A)]^{<\omega}, i \in \{0,1\}\}
$$
$$
\cup \{\{(x,i)\} : x \in \pi(A), i \in \{0,1\}\},
$$
which is obviously countable. We claim that the both assignments $A \to \hat{r}_A$ and $A \to \hat{\mathcal{N}}(A)$ witness that $AD(X)$ is a monotonically retractable space. In fact, it is standard to verify that $\hat{\mathcal{N}}$ is $\omega$-monotone. To prove that $\hat{\mathcal{N}}(A)$ is a network of $\hat{r}_A$ we fix $(x,i) \in AD(X)$ and assume that $\hat{r}_A((x,i)) \in V$ for some open set $V \subseteq AD(X)$. We consider two cases.\medskip

\textit{Case 1.}  $x \in \pi(A)$. In this case, we have that $\{(x,i)\}$ satisfies that $\{(x,i)\} \in \hat{\mathcal{N}}(A)$ and $\hat{r}_A(\{(x,i)\}) = \{(x,i)\} \subseteq V$.\medskip

\textit{Case 2.}  $x \in X \setminus \pi(A)$. We then have that  $\hat{r}_A((x,i)) = (r_{\pi(A)}(x),0) \in V$. Without lost of generality, we may assume that $V = (U \times \{0,1\}) \setminus \{(r_{\pi(A)}(x),1)\}$ where $U$ is a neighborhood of $r_{\pi(A)}(x)$ in $X$. Since $\mathcal{N}(\pi(A))$ is a network of $r_{\pi(A)}$, there is $N \in \mathcal{N}(\pi(A))$ such that $x \in N$ and $r_{\pi(A)}(N) \subseteq U$. If $r_{\pi(A)}(x) \in \pi(A)$, then we put $\hat{N} = (N \setminus \{r_{\pi(A)}(x)\}) \times \{i\}$, and note that $(x,i) \in \hat{N} \in \hat{\mathcal{N}}(A)$ and
\begin{equation*}
\begin{split}
\hat{r}_A(\hat{N}) & \subseteq \big(r_{\pi(A)}(N) \times \{0\}\big) \cup \big(r_{\pi(A)}(N \cap \pi(A)) \times \{1\}\big) \\
& = \big(r_{\pi(A)}(N) \times \{0\}\big) \cup \big((N \cap \pi(A)) \times \{1\}\big) \\
& \subseteq \big(U \times \{0\}\big) \cup \big((U \setminus \{r_{\pi(A)}(x)\}) \times \{1\} \big) \\
& = (U \times \{0,1\}) \setminus \{(r_{\pi(A)}(x),1)\} \\
& = V.
\end{split}
\end{equation*}
For the case $r_{\pi(A)}(x) \not\in \pi(A)$ we set $\hat{N} = N \times \{i\}$, and as above we have that $(x,i) \in \hat{N} \in \hat{\mathcal{N}}(A)$ and $\hat{r}_A(\hat{N}) \subseteq V$.
\endproof

\begin{corollary}\label{CSPB}
If $X$ is the $\Sigma$-product of $\omega_1$-many copies of $\R$, then $C_p(AD(X))$ is Lindel\"of.
\end{corollary}

\proof
 We know from \cite[Corollary 3.14]{rjs1} that $X$ is monotonically retractable. Hence, by applying Theorem \ref{ADMR}, we obtain that  $AD(X)$ is also monotonically retractable. Finally, by \cite[Theorem 3.18]{rjs1}, the space $C_p(AD(X))$ is Lindel\"of.
\endproof

It is proved in \cite{kld} that the Alexandroff duplicate of an Eberlein (a Corson) compact is an Eberlein (a Corson) compact. By a slight modification of the proof of Theorem 2.13 (ii) from \cite{kld}, we can show the following result.

\medskip

\begin{proposition}\label{PADGC}
If $X$ has a weakly $\sigma$-point finite $T_0$-separating family\footnote{Recall that a family $\mathcal{A}$ of subsets of a space $X$ is \textit{$T_0$-separating} if for any distinct points $x,y \in X$ there exists $A \in \mathcal{A}$ such that $A \cap \{x,y\}$ is a singleton. A family $\mathcal{U} = \bigcup_{n \in \mathbb{N}}\mathcal{U}_n$ of subsets of $X$ is called \textit{weakly $\sigma$-point-finite} if  $\mathcal{U} = \bigcup\{\mathcal{U}_n : n \in \mathbb{N}, \mathcal{U}_n \textnormal{ is point finite at } x\}$ for all $x \in X$.} of cozero subsets of $X$, then $AD(X)$ has the same property.
\end{proposition}

{\it Gul'ko compact spaces} are precisely compact spaces for which $C_p(X)$ has the Lindel\"of $\Sigma$-property. G. A. Sokolov proved in \cite{sok}  that a compact space $X$ is Gul'ko compact if and only if $X$ has a weakly $\sigma$-point finite $T_0$-separating family of cozero subsets of $X$. Hence, we have the next corollaries.

\begin{corollary}\label{ADGC}
If $X$ is Gul'ko compact, then $AD(X)$ is also Gul'ko compact.
\end{corollary}

\begin{corollary}\label{CPADGC}
If $X$ is compact and $C_p(X)$ has the Lindel\"of $\Sigma$-property, then $C_p(AD(X))$ has the Lindel\"of $\Sigma$-property.
\end{corollary}

We end this section with two open questions.

\medskip

We have proved that $C_p(AD(X))$ is Lindel\"of whenever $X$ is a Corson compact space. A more general class of compact spaces for which $C_p(X)$ is Lindel\"of is the class of Sokolov spaces with $t(X^n) \leq \omega$ for each $n \in \N$ (for the definition of  Sokolov space see \cite{tka}). This suggests the next natural question.

\begin{question}
Let $X$ be a Sokolov (compact) space such that $t(X^n) \leq \omega$ for each $n \in \N$. It is true that $C_p(AD(X))$ is Lindel\"of.
\end{question}

\begin{question}
Assume that $X$ is pseudocompact and $C_p(X)$ is Lindel\"of $\Sigma$. Is it true that $C_p(AD(X))$ is also Lindel\"of $\Sigma$?
\end{question}

\section{$q$-skeletons and Corson compact spaces}

To start  this section we state the notion of a $q$-skeleton which will be very useful to find Corson compact spaces inside of spaces of continuous functions.

\begin{definition}\label{DRQS} Let $X$ be a space.
A \textit{$q$-skeleton} on  $X$ is a family of $\R$-quotient maps $\{q_s : X \to X_s \mid s \in \Gamma\}$ indexed by an up-directed,  $\sigma$-complete partially ordered nonempty set $\Gamma$ together with a family $\{D_s : s \in \Gamma\}$ of countable subsets of $X$ such that:
\begin{enumerate}[(i)]
	\item the set $q_s(D_s)$ is dense in $X_s$,
  \item if $s,t \in \Gamma$ and $s \leq t$, then there exists a continuous surjective  map $p_{t,s} : X_t \to X_s$ such that $q_s = p_{t,s}\circ q_t$, and
  \item the assignment $s \to D_s$ is $\omega$-monotone.
\end{enumerate}
In addition, if $C_p(X) = \bigcup_{s \in \Gamma}q_s^\ast(C_p(X_s))$, then we say that the $q$-skeleton  is  {\it full}.
\end{definition}

Our task is to prove that countably compact subspaces of $C_p(X)$ have a full $r$-skeleton whenever $X$ has a full $q$-skeleton.
Before, we prove a very technical lemma.

\begin{lemma}\label{LAMN}
Let $K$ be an infinite subspace of $X$ and let $\Gamma$ be an up-directed and $\sigma$-complete partially ordered set. Assume that for each $s \in \Gamma$ we have assigned $\mathcal{M}(s) \in [C_p(X)]^{\leq\omega}$ and $\mathcal{N}(s) \in[K]^{\leq\omega}$ so that:
\begin{enumerate}[(a)]
  \item $\mathcal{M}$ and $\mathcal{N}$ are $\omega$-monotone,
		\item $K = \bigcup_{s \in \Gamma}\mathcal{N}(s) = \bigcup_{s \in \Gamma}\textnormal{cl}_X(\mathcal{N}(s))$, and
	\item $\Delta_{\mathcal{M}(s)} : X \to \R^{\mathcal{M}(s)}$ maps $\textnormal{cl}(\mathcal{N}(s))$ homeomorphically to $\Delta_{ \mathcal{M}(s)}(K)$.
\end{enumerate}
Then $K$ has a full $r$-skeleton.
\end{lemma}
\proof

For each $s \in \Gamma$ we define $r_s = (\Delta_{\mathcal{M}(s)} \restriction_{ \textnormal{cl}(\mathcal{N}(s))})^{-1}\circ (\Delta_{\mathcal{M}(s)} \restriction_K)$. Then $r_s$ is a continuous retraction from $K$ to  $\textnormal{cl}(\mathcal{N}(s))$. We will verify that $\{r_s : s \in \Gamma\}$ satisfies the conditions of  a full $r$-skeleton in $K$.

$(i)$ It is clear that $r_s(X)$ is cosmic for each $s \in \Gamma$.
	
$(ii)$ Assume that $s \leq t$ and $x \in K$. Let $y = r_t(x)$. It follows from $r_t(x) = y = r_t(y)$ that $\Delta_{\mathcal{M}(t)}(x) = \Delta_{\mathcal{M}(t)}(y)$. In particular, we have that  $\Delta_{\mathcal{M}(s)}(x) = \Delta_{\mathcal{M}(s)}(y)$. Hence, $r_s(x) = r_s(y) = r_s(r_t(x)) = r_s \circ r_t(x)$. On the other hand, we know that $r_s(x) \in \textnormal{cl}(\mathcal{N}(s)) \subseteq \textnormal{cl}(\mathcal{N}(t))$ which implies that $r_t \circ r_s(x) = r_t(r_s(x)) = r_s(x)$.
	
$(iii)$ Assume that $\{s_n : n \in \mathbb{N}\} \subseteq \Gamma$ and $s_n \leq s_{n+1}$ for each $n \in \mathbb{N}$. Since $\Gamma$ is $\sigma$-complete, $t = \sup\{s_n\}_{n \in \mathbb{N}}$ exists and $\mathcal{M}(t) = \bigcup\{\mathcal{M}(s_n) : n \in \mathbb{N}\}$. Pick $x \in K$ and let $y = r_t(x)$. Then $\Delta_{\mathcal{M}(t)}(x) = \Delta_{\mathcal{M}(t)}(y)$ and, in particular, we have that  $\Delta_{\mathcal{M}(s_n)}(x) = \Delta_{\mathcal{M}(s_n)}(y)$ for each $n \in \mathbb{N}$. Observe that $\Delta_{\mathcal{M}(s_n)}(r_{s_n}(x)) = \Delta_{\mathcal{M}(s_n)}(x)$ for every  $n \in \mathbb{N}$. It follows that $\Delta_{\mathcal{M}(s_n)}(y) = \Delta_{\mathcal{M}(s_n)}(r_{s_n}(x))$ for each $n \in \mathbb{N}$. Then, by pointwise convergence, we obtain that  $\Delta_{\mathcal{M}(t)}(y) = \lim_{n \to \infty} \Delta_{\mathcal{M}(t)}(r_{s_n}(x))$ and hence $$r_t(x) = r_t(y)  = \lim_{n \to \infty} r_t(r_{s_n}(x)) = \lim_{n \to \infty} r_{s_n}(x).$$
	
$(iv)$ If $x \in \mathcal{N}(s)$ for some $s \in \Gamma$, then $r_s(x) = x$ and so $x = \lim_{s \in \Gamma} r_s(x)$. Finally,  $K = \bigcup_{s \in \Gamma}\textnormal{cl}(\mathcal{N}(s)) = \bigcup_{s \in \Gamma} r_s(X)$.
\endproof

Now we are ready to prove our main result about $q$-skeletons.

\begin{theorem}\label{TCCER}
Assume that the families $\{q_s : s \in \Gamma\}$ and $\{D_s : s \in \Gamma\}$ form a $q$-skeleton on $X$. If a set $K \subseteq \bigcup_{s \in \Gamma}q_s^\ast(C_p(X_s)) \subseteq C_p(X)$ is countably compact, then $K$ has a full $r$-skeleton.
\end{theorem}

\proof
Our strategy is to define two functions $\mathcal{M} : \Gamma \to [C_p(C_p(X))]^{\leq\omega}$ and $\mathcal{N} : \Gamma \to [C_p(X)]^{\leq\omega}$ as in the Lemma \ref{LAMN}. To bring this about, we need to construct tree auxiliary $\omega$-monotone assignments $\mathcal{D} : [C_p(X)]^{\leq\omega} \to [X]^{\leq\omega}$, $\mathcal{E} : [X]^{\leq\omega} \to [C_p(C_p(X))]^{\leq\omega}$ and $\mathcal{A} : [C_p(C_p(X))]^{\leq\omega} \to [C_p(X)]^{\leq\omega}$.\medskip

For each $f \in K$ fix $s_f \in \Gamma$ such that $f \in q_{s_f}^\ast\big(C_p(X_{s_f})\big)$. Consider the function $\gamma : [K]^{\leq\omega} \to \Gamma$ associated to the assignment $f \to s_F$ as in Lemma \ref{LAMG}. If $A \in [K]^{\leq\omega}$ put $\mathcal{D}(A) = D_{\gamma(A)} \in [X]^{\leq\omega}$. It is easy to verify that $\mathcal{D}$ is $\omega$-monotone.\medskip

\textbf{Claim 1.} For each $A \in [K]^{\leq\omega}$ we have $A \subseteq q_{\gamma(A)}^\ast(C_p(X_{\gamma(A)}))$.\medskip

\textit{Proof of Claim 1.} Pick $f \in A$. By the election of $\gamma$, we know that $\gamma(A) \geq \gamma(\{f\}) \geq s_f$. Now, condition (ii) from definition \ref{DRQS} asserts that there exists a continuous surjective function $p_{\gamma(A),s_f} : X_{\gamma(A)} \to X_{s_f}$ such that $q_{s_f} = p_{\gamma(A),s_f}\circ q_{\gamma(A)}$. Hence, $f \in q_{s_f}^\ast(C_p(X_{s_f}))  \subseteq q_{\gamma(A)}^\ast(C_p(X_{\gamma(A)}))$. \medskip

Consider the $\omega$-monotone assignments
$$
\mathcal{W}_0 : [X]^{\leq\omega} \to [\mathcal{P}(C_p(X))]^{\leq\omega} \ \text{and} \ \mathcal{W} : [\mathcal{P}(C_p(X))]^{\leq\omega} \to [\mathcal{P}(C_p(C_p(X)))]^{\leq\omega}
$$
defined as in example \ref{EOW}.  For each nonempty set $N \in \bigcup\{\mathcal{W}(\mathcal{W}_0(D)) : D \in [X]^{\leq\omega}\}$ we fix $\phi_N \in N\subseteq C_p(C_p(X))$. Then, for each  $D \in [X]^{\leq\omega}$ define
$$
\mathcal{E}(D)  = \{\phi_N : N \in \mathcal{W}(\mathcal{W}_0(D)) \ \text{and} \ N \not= \emptyset\} \in [C_p(C_p(X))]^{\leq\omega}.
$$
It is evident that  $\mathcal{E}$ is $\omega$-monotone.
Towards the  definition of the assignment $\mathcal{A}$, for every $G \in [C_p(C_p(X))]^{<\omega}$ fix $A_G \in [K]^{\leq\omega}$ so that $\Delta_G(A_G)$ is dense in $\Delta_G(K)$. For each $E \in [C_p(C_p(X))]^{\leq\omega}$ we define $\mathcal{A}(E) = \bigcup\{A_G : G \in [E]^{\leq \omega}\} \in [K]^{\leq\omega}$. Note that $\mathcal{A}$ is $\omega$-monotone and $\Delta_E(\mathcal{A}(E))$ is dense in $\Delta_E(K)$ for all $E \in [C_p(C_p(X))]^{\leq\omega}$.

For each $A \in [K]^{\leq\omega}$, we define $\mathcal{M}(A) = \mathcal{E}(\mathcal{D}(A)) \in [C_p(C_p(X))]^{\leq\omega}$ and $\mathcal{N}(A) = \mathcal{A}(\mathcal{E}(\mathcal{D}(A))) \in [K]^{\leq\omega}$. It is straightforward to verify that $\mathcal{M}$ and $\mathcal{N}$ are $\omega$-monotone. Let $\overline{\mathcal{N}}(A)$ be the closure of $A$ under $\mathcal{N}$ (see Definition \ref{DCOM}). We know that $\overline{\mathcal{N}}$ is a $\omega$-monotone map. As a consequence the set $\Gamma^\prime = \overline{\mathcal{N}}([K]^{\leq \omega})$ is up-directed and $\sigma$-complete. According to Lemma \ref{LCOM}, we have $A \subseteq \overline{\mathcal{N}}(A) = \mathcal{N}(\overline{\mathcal{N}}(A))$ for all $A \in [K]^{\leq \omega}$ and hence $K = \bigcup\{\overline{\mathcal{N}}(A) : A \in [K]^{\leq \omega}\} = \bigcup_{A \in \Gamma^\prime}\mathcal{N}(A)$. In order to finish the proof, in virtue of Lemma \ref{LAMN} it is enough to show the following two claims.\medskip

\textbf{Claim 2.}  The set $\textnormal{cl}_K(\mathcal{N}(A)) \subseteq q_{\gamma(A)}^\ast(C_p(X_{\gamma(A)}))$ and $\textnormal{cl}_K(\mathcal{N}(A))$ is compact for every $A \in \Gamma^\prime$.\medskip

\textit{Proof of Claim 2}. Observe from Claim $1$ that $\mathcal{N}(A) = A \subseteq K \cap q_{\gamma(A)}^\ast(C_p(X_{\gamma(A)}))$. Since $q_{\gamma(A)}$ is an $\R$-quotient map, the set $q_{\gamma(A)}^\ast(C_p(X_{\gamma(A)}))$ is closed in $C_p(X)$ (for a proof of this fact see \cite[Theorem 0.4.10]{ark}). It follows that $K \cap q_{\gamma(A)}^\ast(C_p(X_{\gamma(A)}))$ is countably compact. On the other hand, since the space $X_{\gamma(A)}$ is separable,  we can apply \cite[Theorem I.1.4]{ark} to see that $C_p(X_{\gamma(A)})$ (and hence $q_{\gamma(A)}^\ast(C_p(X_{\gamma(A)}))$) admits a condensation onto a second countable space. In particular, we have that  $K \cap q_{\gamma(A)}^\ast(C_p(X_{\gamma(A)}))$ admits a condensation onto a second countable space. Since any condensation from a countably compact space onto a Frechet-Urysohn space is a homeomorphism, the space $K \cap q_{\gamma(A)}^\ast(C_p(X_{\gamma(A)}))$ is second countable and hence compact. Therefore, $\textnormal{cl}_K(\mathcal{N}(A)) \subseteq q_{\gamma(A)}^\ast(C_p(X_{\gamma(A)}))$ and $\textnormal{cl}_K(\mathcal{N}(A))$ is compact.\medskip

\textbf{Claim 3.} For each $A \in \Gamma^\prime$ the function $\Delta_{\mathcal{M}(A)}$ maps $\textnormal{cl}_K(\mathcal{N}(A))$ homeomorphically onto $\Delta_{\mathcal{M}(A)}(K)$.\medskip

\textit{Proof of Claim 3}. Since $\Delta_{\mathcal{E}(\mathcal{D}(A))}(\mathcal{A}(\mathcal{E}(\mathcal{D}(A))))$ is dense in $\Delta_{\mathcal{E}(\mathcal{D}(A))}(K)$ and  $\textnormal{cl}_K(\mathcal{N}(A))$ is compact, we have that  $\Delta_{\mathcal{M}(A)}(\textnormal{cl}_K(\mathcal{N}(A))) = \Delta_{\mathcal{M}(A)}(K)$.  To finish the proof of Claim 2, we shall prove that $\Delta_{\mathcal{M}(A)} \restriction \textnormal{cl}_K(\mathcal{N}(A))$ is one-to-one.
Pick two distinct maps $f,g \in \textnormal{cl}_K(\mathcal{N}(A)) \subseteq q_{\gamma(A)}^\ast(C_p(X_{\gamma(A)}))$. Choose $f^\prime,g^\prime \in C_p(X_{\gamma(A)})$ such that $f = f^\prime \circ q_{\gamma(A)}$ and $g = g^\prime \circ q_{\gamma(A)}$. It is clear that $f^\prime \not= g^\prime$. Since $q_{\gamma(A)}(D_{\gamma(A)})$ is dense in $X_{\gamma(A)}$ we can find $x \in D_{\gamma(A)}$ such that $f(x) = f^\prime \circ q_{\gamma(A)}(x) \not= g^\prime \circ q_{\gamma(A)}(x) = g(x)$. Choose disjoint sets $B_1,B_2 \in \mathcal{B}(\R)$ such that  $f(x) \in B_1$ and $g(x) \in B_2$. Consider the sets  $N_1 = [x;B_1]$, $N_2 = [x;B_2]$ and $N = [N_1,N_2;B_1,B_2]$. Observe that $f \in N_1$ and $g \in N_2$. Then $N_1,N_2 \in \mathcal{W}_0(D_{\gamma(A)})$ and $N \in \mathcal{W}(\mathcal{W}_0(D_{\gamma(A)}))$. The set $N$ is nonempty since $\pi_x \in N$. In this way, $\phi_N \in N$ and $\phi_N \in  \mathcal{E}(D_{\gamma(A)}) = \mathcal{E}(\mathcal{D}(A)) = \mathcal{M}(A)$. It follows that $\phi_N(f) \not = \phi_N(g)$ and thus $\Delta_{\mathcal{M}(A)}(f) \not = \Delta_{\mathcal{M}(A)}(g)$.
\endproof

Now, we give a topological condition that implies Corson property on compact subspaces of $C_p(X)$, such result is analogous to Corollary 5.4 in \cite{ban} by replacing the $\Omega$-property by a full $q$-skeleton.

\begin{corollary}\label{CCCP}
If $X$ has a full $q$-skeleton, then any compact subspace of $C_p(X)$ is Corson.
\end{corollary}

\proof
Let $K$ be a compact subspace of $C_p(X)$. It follows from Theorem \ref{TCCER} that $K$ has a full $r$-skeleton. Now, we apply Theorem 3.11 from \cite{mrk2} to obtain that $K$ is a Corson compact space.
\endproof

\begin{corollary}\label{CIC}
Let $K$ be a compact space. If  $C_p(K)$ is a continuous image of $X$ and  $X$ has a full $q$-skeleton, then $K$ is Corson compact.
\end{corollary}

\proof Let $\phi : X \to C_p(K)$ be a continuous surjection and $\psi : K \to C_p(C_p(K))$ be the natural embedding of $K$ into $C_p(C_p(K))$. Since the map $\phi^\ast C_p(C_p(K)) \to C_p(X)$ is an embedding,  the map $\phi^\ast \circ \psi$ is also an embedding. Hence, we can assume that $K \subseteq C_p(X)$. In virtue of Cororllary \ref{CCCP}, we conclude that $K$ is Corson compact.
\endproof

Now we shall proof that many topological spaces have a full $q$-skeleton. Before we prove a technical lemma.

\begin{lemma}\label{LDRQ}
If  $A \subseteq C_p(X)$, then $A \subseteq \Delta_A^\ast(C_p(\Delta_A(X)))$. If in addition $A$ is a  dense subset of $\Delta_A^\ast(C_p(\Delta_A(X)))$, then the map $\Delta_{\textnormal{cl}(A)}$ is $\R$-quotient.
\end{lemma}

\proof For the first assertion, pick $f \in A$ and let $p_f$ be the projection of $\Delta_A(X)$ on the $f$-th coordinate. Then, we obtain that  $f = p_f \circ \Delta_A \in \Delta_A^\ast(C_p(\Delta_A(X)))$.  For the second part observe that for distinct  points $x,y \in X$ we have $\Delta_{\textnormal{cl}(A)}(x) = \Delta_{\textnormal{cl}(A)}(y)$ iff $\Delta_A(x) = \Delta_A(y)$. As a consequence the natural projection $p_A : \Delta_{\textnormal{cl}(A)}(X) \to \Delta_A(X)$ is a condensation. It follows that $p_A^\ast(C_p(\Delta_A(X)))$ is dense in $C_p(\Delta_{\textnormal{cl}(A)}(X))$. In addition, we also know that $\Delta_A = p_A \circ \Delta_{\textnormal{cl}(A)}$. So we have that
$$
\Delta_A^\ast(C_p(\Delta_A(X))) = \Delta_{\textnormal{cl}(A)}^\ast(p_A^\ast(C_p(\Delta_A(X)))) \subseteq \Delta_{\textnormal{cl}(A)}^\ast(C_p(\Delta_{\textnormal{cl}(A)}(X))).
$$
As a consequence the set $\Delta_A^\ast(C_p(\Delta_A(X)))$ is dense in $\Delta_{\textnormal{cl}(A)}^\ast(C_p(\Delta_{\textnormal{cl}(A)}(X)))$. By our hypothesis, we conclude that $A$ is dense in $\Delta_{\textnormal{cl}(A)}^\ast(C_p(\Delta_{\textnormal{cl}(A)}(X)))$. Therefore, we have that $\Delta_{\textnormal{cl}(A)}^\ast(C_p(\Delta_{\textnormal{cl}(A)}(X))) = \textnormal{cl}(A)$.
The map $\Delta_{\textnormal{cl}(A)}$ is $\R$-quotient because of Theorem 0.4.10 from \cite{ark}.
\endproof

\begin{proposition}\label{MWSFQS}
Every monotonically $\omega$-stable space has a full $q$-skeleton.
\end{proposition}

\proof
Suppose that the function $\mathcal{N} : [C_p(X)]^{\leq\omega} \to [\mathcal{P}(X)]^{\leq\omega}$ witness that $X$ is monotonically $\omega$-stable.

For each $F \in [C_p(X)]^{<\omega}$ we  select a countable dense subset $A_F$ of $\Delta_{F}^\ast(C_p(\Delta_F(X)))$ and then  for every $A \in [C_p(X)]^{\leq\omega}$ we let $\mathcal{A}(A) = \bigcup\{A_F : F \in [A]^{\leq\omega}\}$.\medskip

\textbf{Claim.} $\mathcal{A}(A)$ is a dense subset of $\Delta_{A}^\ast(C_p(\Delta_A(X)))$  for all $A \in [C_p(X)]^{\leq\omega}$.

\medskip

Fix $A \in [C_p(X)]^{\leq\omega}$.  To prove the claim first note that
$$
\mathcal{A}(A) = \bigcup\{A_F : F \in [A]^{\leq\omega}\} \subseteq \bigcup\{\Delta_{F}^\ast(C_p(\Delta_F(X))) : F \in [A]^{\leq\omega}\} \subseteq \Delta_{A}^\ast(C_p(\Delta_A(X))).
$$
Now set a canonical open set $U = [x_1,\ldots,x_n;B_1,\ldots,B_n]$ in $C_p(X)$ and $f \in U \cap \Delta_{A}^\ast(C_p(\Delta_A(X)))$. Then there exists $g \in C_p(\Delta_A(X))$ such that $f = g \circ \Delta_A$. Choose $F \in [A]^{\leq\omega}$ so that $\Delta_F(x_i) = \Delta_F(x_j)$ iff $\Delta_A(x_i) = \Delta_A(x_j)$, for $i,j = 1,\ldots,n$. Select $h \in C_p(\Delta_F(X))$ such that $h(\Delta_F(x_i)) = g(\Delta_A(x_i)) = f(x_i)$ for $i = 1,\ldots,n$. Then $h \circ \Delta_F \in U \cap \Delta_{F}^\ast(C_p(\Delta_F(X)))$. Since the set $A_F$ is dense in $\Delta_{F}^\ast(C_p(\Delta_F(X)))$ we must have that $U \cap \mathcal{A}(A) \supseteq U \cap A_F \not= \emptyset$. Therefore, $\mathcal{A}(A)$ is a dense subset of $\Delta_{A}^\ast(C_p(\Delta_A(X)))$.

\medskip

Consider the closure $\overline{\mathcal{A}}(A)$ of $A$ under $\mathcal{A}$ which is $\omega$-monotone because of Lemma \ref{LCOM}. On the other hand,  for each nonempty set $N \in \bigcup\{\mathcal{N}(A) : A \in [C_p(X)]^{\leq\omega}\}$ fix a point $x_N \in N$. For each $A \in [C_p(X)]^{\leq\omega}$ we define $\mathcal{D}(A) = \{x_N : N \in \mathcal{N}(A) \setminus \{\emptyset\}\}$. It is clear that the map $\mathcal{D}$ is $\omega$-monotone. Consider the up-directed and $\sigma$-compete partially ordered set $\Gamma = \overline{\mathcal{A}}([C_p(X)]^{\leq\omega})$. For each $A \in \Gamma$ let $D_A := \mathcal{D}(A)$ and $q_A := \Delta_{\textnormal{cl}(A)}$.
Fix $A \in \Gamma$. Since $\mathcal{N}(A)$ is a network for $q_A$, the set $q_A(D_A)$ is dense in $q_A(X)$. By the Claim we have that $A = \mathcal{A}(A)$ is dense in $\Delta_{A}^\ast(C_p(\Delta_A(X)))$. According to  Lemma \ref{LDRQ}   $q_A$ is an $\R$-quotient map. For each $B \in \Gamma$ with $A \subseteq B$ let us consider the natural projection $p_{B,A}$ from $q_B(X)$ onto $q_A(X)$ which satisfies $q_A = p_{B,A} \circ q_B$. Finally it is clear from the first part of Lemma \ref{LDRQ} that $C_p(X) = \bigcup\{A : A \in \Gamma\} = \bigcup\{q_A^\ast(C_p(q_A(X))) : A \in \Gamma\}$. Therefore, $X$ has a full $q$-skeleton.
\endproof

It was proved in \cite[Proposition 4.4]{rjs2} and \cite[Corollary 4.14]{rjs2} that Lindel\"of $\Sigma$-spaces and pseudocompact spaces are monotonically $\omega$-stable. Hence we have the following corollary.

\begin{corollary}
If $X$ is either Lindel\"of $\Sigma$ or pseudocompact, then $X$ has a full $q$-skeleton.
\end{corollary}

Now we give conditions under which a $q$-skeleton must be full.

\begin{proposition}\label{LLRQS}
Assume that the families $\{q_s: X \to X_s | s \in \Gamma\}$ and $\{D_s : s \in \Gamma\}$ form a $q$-skeleton on $X$ such that:
\begin{enumerate}
\item $X = \bigcup\{D_s : s \in \Gamma\}$,
\item $q_s \restriction D_s$ is one-to-one, and
\item $l(X^n) \leq \omega$ for all $n \in \N$.
\end{enumerate}
Then the $q$-skeleton is full.
\end{proposition}

\proof
  It is suffices to prove that $C_p(X) = \bigcup_{s \in \Gamma}q_s^\ast(C_p(X_s))$. Given $F \in [X]^{<\omega}$ fix $s_F \in \Gamma$ such that $F \subseteq D_{s_F}$. For each canonical non-empty open set $W = [x_1,\ldots,x_n;B_1,\ldots,B_n] \in \mathcal{W}(F)$ choose a map $h_{F,W} : X_{s_F} \to \R$ such that $h_{F,W}(q_{s_F}(x_i)) \in B_i$ for $i = 1,\ldots,n$. Consider the map $g_{F,W} = h_{F,W} \circ q_{s_F} \in q_{s_F}^\ast(C_p(X_{s_F})) \cap W$. Define
  $$
  C_F = \{g_{F,W} : W \in \mathcal{W}(F), W \not = \emptyset\} \subseteq q_{s_F}^\ast(C_p(X_{s_F})) \subseteq C_p(X).
  $$
Next, for every $A \in [X]^{\leq\omega}$ let $\mathcal{C}(A) = \bigcup\{C_F : F \in [A]^{\leq \omega}\} \in [C_p(X)]^{\leq\omega}$. It is evident that $\mathcal{C}$ is $\omega$-monotone. It follows directly from the construction $\mathcal{C}(X) = \bigcup\{\mathcal{C}(A) : A \in [X]^{\leq\omega}\}$ is dense in $C_p(X)$.

Pick $f \in C_p(X)$. Since $\mathcal{C}(X)$ is dense in $C_p(X)$ and $t(X) = \sup\{l(X^n) : n \in \N\} = \omega$, we can find a countable set $C \subseteq \mathcal{C}(X)$ such that $f \in \textnormal{cl}(C)$. Since $\mathcal{C}$ is $\omega$-monotone we can assume that $C \subseteq \mathcal{C}(A)$ for some $A \in [X]^{\leq\omega}$. Let $s$ be an upper bound of $\{s_F: F \in [A]^{\leq\omega}\}$ and note that
$$
\mathcal{C}(A) = \bigcup\{C_F : F \in [A]^{\leq \omega}\} \subseteq \bigcup\{q_{s_F}^\ast(C_p(X_{s_F})) : F \in [A]^{\leq \omega}\} \subseteq q_{s}^\ast(C_p(X_{s})).
$$
Since $q_{s}$ is an $\R$-quotient map, $f \in \textnormal{cl}(C) \subseteq \textnormal{cl}(\mathcal{C}(A)) \subseteq \textnormal{cl}(q_{s}^\ast(C_p(X_{s}))) = q_{s}^\ast(C_p(X_{s}))$.
\endproof

We finish this section with two open questions.

\begin{question}
Assume that $X$ is Lindel\"of and has a full $q$-skeleton. It is true that each continuous image of $X$ has a full $q$-skeleton?
\end{question}

\begin{question}
Assume that $X$ is Lindel\"of and has a full $q$-skeleton. Must $X^n$ be necessarily Lindel\"of, for each $n \in \N$?
\end{question}

\section{Strong $r$-skeletons}

A space $X$ is monotonically $\omega$-stable iff $C_p(X)$ is monotonically $\omega$-monolithic \cite{rjs2}. We also know that a space $X$ is monotonically retractable iff $C_p(X)$ is monotonically Sokolov \cite{rjs-tka}. Besides, it was proved in \cite{cas-rjs} that a space $X$ is monotonically retractable iff it is monotonically $\omega$-stable and has a full $r$-skeleton:
$$
\begin{array}{cccccc}
X      & \text{monotonically retractable} & = & \text{monotonically $\omega$-stable}     & + & \text{full $r$-skeleton}\\
       & \Updownarrow                     &   & \Updownarrow                             &   & \\
C_p(X) & \text{monotonically Sokolov}     & = & \text{monotonically $\omega$-monolithic} & + & ?\\
\end{array}
$$
One of the purposes of this section is to find a  system of retractions which completes this diagram.
In this way, it is natural to conjecture that a space $X$ is monotonically Sokolov iff it is monotonically $\omega$-monolithic and has a full $r$-skeleton. However this conjecture is false according to Example \ref{EJEM} (c). By this reason, it was asked in  \cite[Question 4.6]{cas-rjs} if there exists a characterization of the monotone Sokolov property using monotone $\omega$-monolithicity and a system of retractions. To solve this problem we will introduce the notion of strong $r$-skeleton and show that a space $X$ is monotonically Sokolov iff it is monotonically $\omega$-monolithic and has a strong $r$-skeleton.

\begin{definition}\label{DSRS}
 We say that a $r$-skeleton$\{r_s : s \in \Gamma\}$ in a space $X$ is a \textit{strong $r$-skeleton} if satisfies the following condition: for each $s \in \Gamma$, $n \in \N$ and a closed subset $F$  of $X^n$ there exists $t \in \Gamma$ such that $s \leq t$ and $r_t^n(F) \subseteq F$, where $r_t^n : X^n \to X^n$ is the n-th power of the map $r_t$.
\end{definition}

It is easy to see that any strong $r$-skeleton is a full $r$-skeleton (the Example \ref{EJEM} (c) assetrts that these two notions are distinct). Besides it is also easy to verify that strong $r$-skeletons are preserved by countable disjoint topological unions and inherited by closed subspaces.
Next, we shall prove that a space $X$ has a strong $r$-skeleton iff $C_p(X)$ has one. To have this done we first prove two auxiliary lemmas.

\begin{lemma}\label{LSRS}
Let $\{r_s : s \in \Gamma\}$ be a strong $r$-skeleton in a space $X$. If $s \in \Gamma$ and $\mathcal{F}$ is a countable subset of  $\bigcup_{n \in \N}exp(X^n)$, then there exists $t \in \Gamma$ with $s \leq t$ satisfying that $r_t^n(F) \subseteq F$ whenever $F \in \mathcal{F}$ and $F \subseteq X^n$.
\end{lemma}

\proof
Suppose that $s \in \Gamma$ and $\mathcal{F} \subseteq \bigcup_{n \in \N} exp(X^n)$ is  countable. Let $\{F_m : m \in \mathbb{N}\}$ be a numeration of $\mathcal{F}$ where each element appears infinitely many times. Choose inductively an increasing sequence $\{s_m : m \in \mathbb{N}\} \subseteq \Gamma$ as follows. Let $s_0 = s$. Assume that $s_k \in \Gamma$ has been defined for each $0 \leq k <m$. Since $F_{m}$ is a closed subset of $X^{n}$ for some $n \in \N$, we can find $s_{m} \in \Gamma$ with $s_k \leq s_m$ for each $k < m$ and $r_{s_m}^{n}(F_m) \subseteq F_m$. We claim that $t = \sup\{s_m\}_{m \in \mathbb{N}}$ is the required element. Indeed, given $F \in \mathcal{F}$ let $N_F = \{m \in \N : F_m = F\}$. Since $N_F$ is cofinal in $\N$, $t = \sup \{s_m : m \in N_F\}$. Note that if $m \in N_F$ and $F$ is a subset of $X^n$, then $r_{s_{m}}^n(F) \subseteq F$. Thus, we obtain that $r_t^n(F) \subseteq F$.
\endproof

\begin{lemma}\label{LFRS}
Assume that $X$ has a full $r$-skeleton $\{r_s : s \in \Gamma\}$. For each  $s \in \Gamma$, we define
$$
\hat{r}_s = r_s^\ast \circ \pi_{r_s(X)} : C_p(X) \to C_p(X).
$$
Then we have that $\{\hat{r}_s : s \in \Gamma\}$ is a $r$-skeleton in $C_p(X)$.
\end{lemma}

\proof
It is easy to verify that each $\hat{r}_s$ is a continuous retraction. We will verify that $\{\hat{r}_s : s \in \Gamma\}$ satisfies conditions (i)-(iv) from Definition \ref{DRS}.\medskip

$(i)$ Pick $s \in \Gamma$. Since the space $r_s(X)$ is cosmic, it is well known that the space $C_p(r_s(X))$ is also cosmic. Hence,  $\pi_{r_s(X)}(C_p(X))$ must be cosmic. Since the map $r_s^\ast$ is an embedding, the space $\hat{r}_s(X) = r_s^\ast(\pi_{r_s(X)}(C_p(X)))$ is also cosmic.\medskip

$(ii)$ If $s,t \in \Gamma$ and $s \leq t$, then for each $f \in C_p(X)$ we have $\hat{r}_s(f) = f \circ r_s = f \circ r_t \circ r_s = \hat{r}_t(f) \circ r_s = \hat{r}_s \circ \hat{r}_t(f)$ and $\hat{r}_s(f) = f \circ r_s = f \circ r_s \circ r_t = \hat{r}_s(f) \circ r_t = \hat{r}_t \circ \hat{r}_s(f)$. Thus, we obtain that $\hat{r}_s  =  \hat{r}_S \circ \hat{r}_t = \hat{r}_t \circ \hat{r}_s$.\medskip

$(iii)$ Given $\{s_n : n \in \mathbb{N}\} \subseteq \Gamma$ with $s_n \leq s_{n+1}$ for each $n \in \mathbb{N}$, we set $t = \sup\{s_n : n \in \mathbb{N}\}$. For each $f \in C_p(X)$ and $x \in X$ we know that $\hat{r}_t(f)(x) = f(r_t(x)) = f(\lim_{n \to \infty} r_{s_n}(x)) = \lim_{n \to \infty} f(r_{s_n}(x)) = \lim_{n \to \infty} \hat{r}_{s_n}(f)(x)$. It follows that $\hat{r}_t(f) = \lim_{n \to \infty} \hat{r}_{s_n}(f)$.\medskip

$(iv)$ Assume that $f \in C_p(X)$ and $f \in U = [x_1,\ldots,x_n;B_1,\ldots,B_n]$ for some canonical open subset of $C_p(X)$. Choose $s \in \Gamma$ such that $r_s(x_i) = x_i$ for $i = 1,\ldots,n$. Assume that $t \in \Gamma$ and $s \leq t$. Then $\hat{r}_t(f)(x_i) = f \circ r_t(x_i) = f(x_i) \in B_i$, for $i = 1,\ldots,n$. It follows that $\hat{r}_t(f) \in U$. Therefore, $f = \lim_{s \in \Gamma}\hat{r}_s(f)$.
\endproof

\begin{theorem}\label{SRSCP}
A space $X$ has a strong $r$-skeleton iff $C_p(X)$ has a strong $r$-skeleton.
\end{theorem}

\proof We first prove the necessity.
Assume that $\{r_s : s \in \Gamma\}$ is a strong $r$-skeleton on $X$. We know  that $\{r_s : s \in \Gamma\}$ is a full $r$-skeleton in $X$. If $\{\hat{r}_s : s \in \Gamma\}$ is as in Lemma \ref{LFRS}, then $\{\hat{r}_s : s \in \Gamma\}$ is an $r$-skeleton in $C_p(X)$. In order to prove that $\{\hat{r}_s : s \in \Gamma\}$ is a strong $r$-skeleton, it is enough to show the following claim.\medskip

\textbf{Claim.} If $s \in \Gamma$, $n \in \N$ and $G$ is a closed subset of $C_p(X)^n$, then there exists $t \in \Gamma$ such that $s \leq t$ and $r_t^n(G) \subseteq G$.\medskip

Let $nX$ be the disjoint topological union of $n$ copies of $X$. For each $s \in \Gamma$ let $nr_{s}$ denote the natural continuous retraction induced by $r_s$ on $nX$. It is standard to verify that $\{nr_s : s \in \Gamma\}$ is a strong $r$-skeleton in $nX$. We identify $C_p(X)^n$ with $C_p(nX)$. For each $m \in \N$ and $B_1,\ldots,B_m \in \mathcal{B}(\R)$ we define
$$
F(G,B_1,\ldots,B_m) = \{(x_1,\ldots,x_m) \in (nX)^m : [x_1,\ldots,x_m;B_1,\ldots,B_m] \cap G = \emptyset\}.
$$
It is easy to verify that $F(G,B_1,\ldots,B_m)$ has an open complement in $(nX)^m$, and hence, it is a closed subset of $(nX)^m$. By Lemma \ref{LSRS} we can find $t \in \Gamma$ such that $s \leq t$ and $(nr_t)^m(F(G,B_1,\ldots,B_m)) \subseteq F(G,B_1,\ldots,B_m)$ for each $m \in \N$ and $B_1,\ldots,B_m \in \mathcal{B}(\R)$. We assert that $t$ is as promised. Assume on the contrary that there exists $f \in G$ such that $\hat{r}_t^n(f) \not \in G$. We identify $\hat{r}_t^n$ with the map $\widehat{nr}_t= (nr_t)^\ast \circ \pi_{nr_t(nX)} : C_p(nX) \to C_p(nX)$. Since $\widehat{nr}_t(f) \not \in G$, then we can choose $x_1,\ldots,x_m \in nX$ and $B_1,\ldots,B_m \in \mathcal{B}(\R)$ such that $\widehat{nr}_t(f) \in [x_1,\ldots,x_m;B_1,\ldots,B_m]$ and $[x_1,\ldots,x_m;B_1,\ldots,B_m] \cap G = \emptyset$. Then we obtain that  $f \circ nr_t \in [x_1,\ldots,x_m;B_1,\ldots,B_m]$ and $(x_1,\ldots,x_m) \in F(G,B_1,\ldots,B_m)$. It follows that $f \in [nr_t(x_1),\ldots,nr_t(x_m);B_1,\ldots,B_m]$. By the properties of $nr_t$ we have that $(nr_t(x_1),\ldots,nr_t(x_m)) \in F(G,B_1,\ldots,B_m)$; that is,
$$
[nr_t(x_1),\ldots,nr_t(x_m);B_1,\ldots,B_m] \cap G = \emptyset,
$$
but this is a contradiction since $f \in G$.

\medskip

For the  other implication of the theorem, assume that $C_p(X)$ has a strong $r$-skeleton. By the firs part $C_p(C_p(X))$ has a strong $r$-skeleton.  Since strong $r$-skeletons are inherited by closed subspaces and $X$ is a closed subspace of $C_p(X)$, we conclude that $X$ has a strong $r$-skeleton.
\endproof

No we will prove the main results of this section. To do that we need some technical assertions.

\medskip

The next Lemma is established by using some basic ideas in the proof of Corolary 4.8 from \cite{rjs-tka}. We give a sketch of the proof.

\begin{lemma}\label{LMR}
If the assignments $A \to r_A$ and $A \to \mathcal{N}(A)$ witnesses that the space $X$ is monotonically retractable, then for every $A \in [X]^{\leq\omega}$ and every $F \in exp(X)$ there exists $B \in [X]^{\leq\omega}$ such that $A \subseteq B$ and $r_B(F) \subseteq F$.
\end{lemma}

\proof Let $A \in [X]^{\leq\omega}$ and  $F \in exp(X)$.
By induction we construct a sequence of countable sets $\{B_n : n \in \N\}$ as follows: $B_0 = A$ and suppose that $B_n$ has been constructed for some $n \in \N$. Fix a countable set $A_n \subseteq X$ such that $A_n \cap F \cap N \not = \emptyset$ whenever $F \cap N \not = \emptyset$, for each $N \in \mathcal{N}(B_n)$. Define $B_{n+1} = A_n \cup B_n$. Finally it can be verified that $B = \bigcup\{B_n : n \in \N\}$ is as promised.
\endproof

\begin{theorem}\label{MRSRS}
If $X$ is monotonically retractable space, then $X$ a strong $r$-skeleton.
\end{theorem}

\proof
Assume that the assignments $A \to r_A$ and $A \to \mathcal{N}(A)$ witnesses that $X$ is monotonically retractable. As in the proof of \cite[Theorem 4.3]{cas-rjs}, for each $A \in [X]^{\leq\omega}$ we can find $E(A) \in [X]^{\leq\omega}$ so that $A \subseteq E(A)$, the assignment $A \to E(A)$ is $\omega$-monotone, and   $\{r_s : s \in \Gamma\}$ is a full $r$-skeleton in $X$, where $\Gamma = \{E(A) : A \in [X]^{\leq\omega}\}$. To see that $\{r_s : s \in \Gamma\}$ is a strong $r$-skeleton, we only need to show the following claim.\medskip

\textbf{Claim.} If $s \in \Gamma$, $n \in \N$ and $F$ is a closed subset of $X^n$, then there exists $t \in \Gamma$ such that $s \leq t$ and $r_t^n(F) \subseteq F$.\medskip

{\it Proof of the Claim.} We know that $X^n$ is monotonically retractable \cite{rjs1}, which can be justified as follows:  For every  countable set $B \subseteq X^n$ we consider
\begin{enumerate}
 \item $S(B) = \bigcup_{i = 1}^{n}p_i(B)$, where $p_i$ is the $i$th-projection,

  \item $\hat r_B = r_{E(S(B))}^n$ and

  \item $\hat{\mathcal{N}}(B) = \{\prod_{i=1}^{n}N_i : N_1,\ldots,N_n \in \mathcal{N}(E(S(B)))\}$.
  \end{enumerate}
  Indeed, the assignments $B \to \hat r_B$ and $B \to \hat{\mathcal{N}}(B)$ witnesses that $X^n$ is monotonically retractable. Choose $A \in [X]^{\leq\omega}$ for which  $s = E(A)$. In virtue of Lemma \ref{LMR}, there exists a countable set $B \subseteq X^n$ such that $A^n \subseteq B$ and $\hat  r_B(F) \subseteq F$. It follows that if $t = E(S(B))$, then $s \leq t$ and $r_t^n(F) \subseteq F$.
\endproof

\begin{corollary}\label{CMRSRS}
A space $X$ is monotonically retractable iff is monotonically $\omega$-stable and has a strong $r$-skeleton.
\end{corollary}

In a dual way, we have the following Corollary.

\begin{corollary}\label{CMSSRS}
A space $X$ is monotonically Sokolov iff is monotonically $\omega$-monolithic and has a strong $r$-skeleton.
\end{corollary}

We know that if $X$ is monotonically Sokolov, then $C_p(X)$ is monotonically retractable and hence any compact subspace of $C_p(X)$ is Corson. In a similar way, it is also known that if $X$ is monotonically retractable, then $C_p(X)$ is monotonically Sokolov and hence any compact subspace of $C_p(X)$ is Corson. So it is natural to ask if any compact subspace of $C_p(X)$ is Corson when $X$ has a strong $r$-skeleton. The following theorem provides a positive answer to this question because of Corollary \ref{CCCP}.

\begin{theorem}\label{SRSFQS}
If $X$ has a strong $r$-skeleton, then $X$ has a full $q$-skeleton.
\end{theorem}

\proof
Let $\{r_s : s \in \Gamma\}$ be a strong $r$-skeleton in $X$. In order get a full $q$-skeleton in $X$ we will construct a map $\gamma : [exp(X)]^{\leq\omega} \to \Gamma$ and a map $\mathcal{C} : [exp(X)]^{\leq\omega} \to [exp(X)]^{\leq\omega}$.

For each $\mathcal{F} \in [exp(X)]^{<\omega}$ we will assign an element $\gamma(\mathcal{F}) \in \Gamma$, by induction on the cardinality of $\mathcal{F}$, as follows: $\gamma(\emptyset) = s_0$ for some arbitrary $s_0 \in \Gamma$. Given a non-empty family $\mathcal{F} \in [exp(X)]^{<\omega}$ choose $\gamma(\mathcal{F}) \in \Gamma$ such that $\gamma(\mathcal{F}^\prime) \leq \gamma(\mathcal{F})$ for each $\mathcal{F}^\prime \subseteq \mathcal{F}$ with $\mathcal{F}^\prime \not= \mathcal{F}$; and $r_{\gamma(\mathcal{F})}(F) \subseteq F$ for each $F \in \mathcal{F}$. Now, for each $\mathcal{F} \in [exp(X)]^{\leq\omega}$ let $\gamma(\mathcal{F}) = \sup\{\gamma(\mathcal{F}^\prime) : \mathcal{F}^\prime \in [\mathcal{F}]^{< \omega}\}$. Note, in the general case,  that $\gamma(\mathcal{F})$ is well defined and
\begin{enumerate}[(a)]
  \item $r_{\gamma(\mathcal{F})}(F) \subseteq F$ for each $F \in \mathcal{F}$;
	\item if $\mathcal{F} \subseteq \mathcal{F}^\prime$, then $\gamma(\mathcal{F}) \leq \gamma(\mathcal{F}^\prime)$; and
	\item if $\{\mathcal{F}_n : n \in \mathbb{N}\} \subseteq [exp(X)]^{\leq\omega}$, $\mathcal{F}_n \subseteq \mathcal{F}_{n+1}$ for each $n \in \mathbb{N}$ and $\mathcal{F} = \bigcup_{n \in \mathbb{N}}\mathcal{F}_n$; then  $\gamma(\mathcal{F}) =  \sup_{n \in \mathbb{N}} \gamma(\mathcal{F}_n)$.
\end{enumerate}
To construct $\mathcal{C}$, if $\mathcal{F} \in [exp(X)]^{<\omega}$, since the space $r_{\gamma(\mathcal{F})}(X)$ is cosmic we can fix a countable family of singletons $\mathcal{C}_\mathcal{F} \subseteq exp(X)$ such that $\bigcup\mathcal{C}_\mathcal{F}$ is dense in $r_{\gamma(F)}(X)$. For every $\mathcal{F} \in [exp(X)]^{\leq\omega}$ let $\mathcal{C}(\mathcal{F}) = \bigcup\{\mathcal{C}_\mathcal{F} : \mathcal{F} \in [\mathcal{F}]^{< \omega}\}$. Consider the closure $\overline{\mathcal{C}}(\mathcal{F})$ of $\mathcal{F}$ under $\mathcal{C}$. Observe that $\mathcal{C}$ and $\overline{\mathcal{C}}$ are $\omega$-monotone.

We are ready to construct a full $q$ skeleton on $X$. The set $\Gamma^\prime = \{\gamma(\overline{\mathcal{C}}(\mathcal{F})) : \mathcal{F} \in [X]^{\leq \omega}\}$ is an up-directed partially ordered subset of $\Gamma$. For every $s = \gamma(\overline{\mathcal{C}}(\mathcal{F})) \in \Gamma^\prime$ we define $q_s = r_s$ and $D_s = \bigcup\mathcal{C}(\overline{\mathcal{C}}(\mathcal{F}))$. The map $q_s$ being a continuous retraction is an $\R$-quotient map. We shall verify that $\{q_s : s \in \Gamma^\prime\}$ and $\{D_s : s \in \Gamma^\prime\}$ form a full $q$-skeleton on $X$. \medskip

$(i)$ We claim that $q_{s}(D_s)$ is dense in $q_s(X)$ whenever $s = \gamma(\overline{\mathcal{C}}(\mathcal{F})) \in \Gamma^\prime$. Assume that $\overline{\mathcal{C}}(\mathcal{F}) = \bigcup\{\mathcal{F}_n : n \in \N\}$ where $\mathcal{F}_n$ is finite and $\mathcal{F}_n \subseteq \mathcal{F}_{n+1}$ for each $n \in \N$. Then $q_s = r_{\gamma(\overline{\mathcal{C}}(\mathcal{F}))} = \lim_{n \to \infty} r_{\gamma(\mathcal{F}_n)}$. As a consequence $q_s(X) = r_{\gamma(\overline{\mathcal{C}}(\mathcal{F}))}(X) \subseteq \textnormal{cl}(\bigcup\{r_{\gamma(\mathcal{F}_n)}(X) : n \in \N\})$. By construction, we have that  $\bigcup\mathcal{C}_{\mathcal{F}_n}$ is dense in $r_{\gamma(\mathcal{F}_n)}(X)$ for each $n \in \N$. It follows that $\bigcup_{n \in \N}\big(\bigcup\mathcal{C}_{\mathcal{F}_n}\big)$ is dense in $\textnormal{cl}(\bigcup\{r_{\gamma(\mathcal{F}_n)}(X) : n \in \N\})$. Since $\bigcup\{\mathcal{C}_{\mathcal{F}_n} : n \in \N\} \subseteq \mathcal{C}(\overline{\mathcal{C}}(\mathcal{F}))$, the set $D_s$ is dense in $\textnormal{cl}(\bigcup\{r_{\gamma(\mathcal{F}_n)}(X) : n \in \N\})$. Hence, $D_s$ is dense in $q_s(X)$ as well. We only need to show that $q_{s}(x) = x$ for each $x \in D_s$. In fact, if $x \in D_s$, then $\{x\} \in \mathcal{C}(\overline{\mathcal{C}}(\mathcal{F})) \subseteq \overline{\mathcal{C}}(\mathcal{F})$. By (a) we obtain that $r_{\gamma(\overline{\mathcal{C}}(\mathcal{F}))}(\{x\}) \subseteq \{x\}$; that is, $q_{s}(x) = x$.\medskip

Conditions $(ii)$ and $(iii)$ in Definition \ref{DRQS} are easy to verify. To finish the proof we must show that the $q$-skeleton is full. In other words, $C_p(X) = \bigcup_{s \in \Gamma^\prime}q_s^\ast(C_p(X_s))$ where $X_s = q_s(X)$ for $s \in \Gamma^\prime$. Let $f \in C_p(X)$ and consider the family $\mathcal{F} = \{f^{-1}(\textnormal{cl}(B)): B \in \mathcal{B}(\R)\}$. If $s = \gamma(\overline{\mathcal{C}}(\mathcal{F})) \in \Gamma^\prime$, then the continuous map $q_s : X \to X$ satisfies $q_s(F) \subseteq F$ for any $F \in \mathcal{F}$. According to Lemma 4.15 of \cite{rjs-tka}, $f = (f \restriction_{X_s}) \circ q_s$, that is, $f \in q_s^\ast(C_p(X_s))$.
\endproof

\begin{corollary}\label{CCPC}
A compact space $X$ is Corson iff $C_p(X)$ has a full $q$-skeleton.
\end{corollary}

\proof
Assume that $X$ is a Corson compact. By \cite[Corollary 4.10]{cas-rjs} and Corollary \ref{CMSSRS} the space $C_p(X)$ has a strong $r$-skeleton. So we can apply Theorem \ref{SRSFQS} to see that $C_p(X)$ has a full $q$-skeleton. The other implication follows from Corollary \ref{CCCP} and the fact that $X$ can be embedded in $C_p(C_p(X))$.
\endproof

The following diagram represents the relationships among the topological properties that so far have been studied.

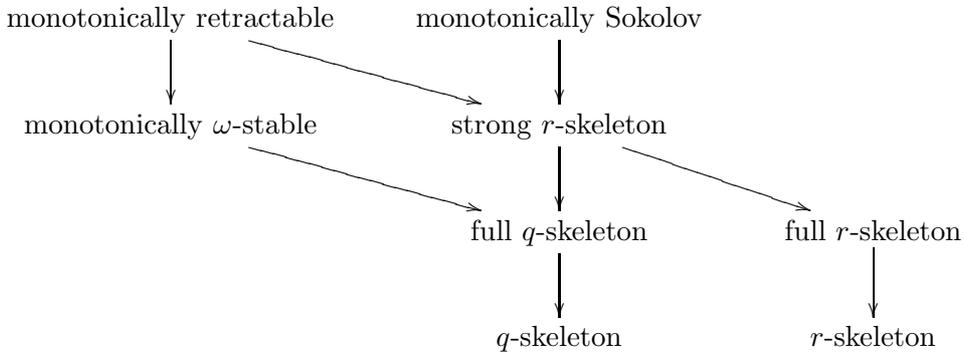
\begin{displaymath}
\xymatrix
{
\textnormal{monotonically retractable}\ar[d]\ar[rd]& \textnormal{monotonically Sokolov}\ar[d]      &\\
\textnormal{monotonically $\omega$-stable}  \ar[rd]& \textnormal{strong $r$-skeleton}\ar[d] \ar[rd]&\\
                                                   & \textnormal{full $q$-skeleton}\ar[d]          &\textnormal{full $r$-skeleton}\ar[d]\\                                                    & \textnormal{$q$-skeleton}                     &\textnormal{$r$-skeleton}
}
\end{displaymath}

Now we list  examples to show that, in general, no one of these implications can be reversed.

\begin{example}\label{EJEM}
\begin{enumerate}[(a)]
\item The compact space $[0,1]^{\omega_1}$ is not Corson, so it is an example of a monotonically $\omega$-stable space which is not monotonically retractable (see \cite[Corollary 4.14]{rjs2} and \cite[Corollary 1.2]{mrk-kld}). According to Proposition \ref{MWSFQS} and \cite[Theorem 3.11]{mrk2}, this space also provides an example of a space with a full $q$-skeleton which does not have  a strong $r$-skeleton.

\item The linearly ordered space of all countable ordinals $\omega_1$ is monotonically retractable hence has a strong $r$-skeleton (see \cite[Theorem 3.8]{rjs1} and Theorem \ref{MRSRS}). However $\omega_1$ is not monotonically Sokolov because it is not Lindel\"of (see \cite[Corollary 4.20]{rjs-tka}). In addition, Theorem \ref{SRSCP} and \cite[Corollary 4.19]{rjs-tka} imply that $C_p(\omega_1)$ has a strong $r$-skeleton but it is not monotonically retractable.

\item Fix $\kappa > \omega$ and let $D_\kappa$ be the set $\kappa$ endowed with the discrete topology. Then $D_\kappa$ has a full $r$-skeleton but does not have a strong $r$-skeleton. To see that $D_\kappa$ has a full $r$-skeleton, consider the up-directed and $\sigma$-complete partially ordered set $\Gamma = [\kappa]^{\leq\omega}$. For each $A \in \Gamma$ we define $m(A) = \min A$ and consider the continuous retraction $r_A : D_\kappa \to D_\kappa$ given by $r_A(\alpha) = \alpha$ if $\alpha \in A$ and $r_A(\alpha) = m(A)$ otherwise. It is easy to verify that $\{r_A : A \in \Gamma\}$ is a full $r$-skeleton in $D_\kappa$. On the other hand, note that $[0,1]^{\kappa} \subseteq \R^{\kappa} = C_p(D_\kappa)$. Since $[0,1]^{\kappa}$ is not a Corson compact, by Corollary \ref{CCCP}, the space $D_\kappa$ does not have a full $q$-skeleton and, as a consequence, it does not admit a strong $r$-skeleton (see Theorem \ref{SRSFQS}).

\item By \cite[Proposition 4.29]{rjs-tka} the space $L_\kappa$ is monotonically Sokolov. Applying  Corollary \ref{CMSSRS} and Theorem \ref{SRSFQS} we obtain that it has a full $q$-skeleton. On the other hand, according to Corollary \ref{CMRSRS}, the space $L_\kappa$ can not be monotonically $\omega$-stable since it does not have countable tightness (see \cite[Corollary 3.20]{rjs1}).

\item A Valdivia compact space which is not Corson is an example of a space with an $r$-skeleton which does not have a full retractional skeleton (see \cite[Theorem 6.1]{kbs-mich} and \cite[Theorem 3.11]{mrk2}).

\item Let $K$ be a compact space which is not Corson. It is easy to construct, by using only one projection, a trivial $q$-skeleton in $C_p(K)$. However, by Corollary \ref{CCPC}, the space $C_p(K)$ does not have a full $q$-skeleton.
\end{enumerate}
\end{example}

As we have mentioned above,  a compact space is Corson iff it has a full $r$-skeleton iff it is monotonically retractable, and a compact space is Valdvia iff it has a commutative $r$-skelton.

\begin{problem}
Find a  characterization of Valdivia compact spaces  similar to monotone retractability.
\end{problem}

\section{$q$-skeletons in $L_\kappa \times K$ for a compact $K$}

In the following, we use $q$-skeletons to prove  that a compact space $X$ is Corson whenever $C_p(X)$ is the continuous image of a closed subspace of $L_\kappa^\omega \times K$, where $K$ is a
compact space. In virtue of Corollary \ref{CIC}, it suffices  to show that each closed subspace of $L_\kappa^\omega \times K$ have a full $q$-skeleton.
First we prove a technical lemma and introduce some notation.

\begin{lemma}\label{LRQM}
Assume that $\phi : K \to K^\prime$ is a continuous surjective map, $r : Y \to Y$ is a continuous retraction, $D \subseteq X \subseteq Y \times K$ and $q = r \times \phi \restriction_X$. If $K$ is compact, $X$ is closed in $Y \times K$, $p_Y(D) \subseteq r(Y)$ and $q(D)$ is dense in $q(X)$, then $q: X \to q(X)$ is an $\R$-quotient map.\medskip
\end{lemma}

\proof
We shall prove that $q$ is the composition of two $\R$-quotient maps. Let $Y^\prime = r(Y)$. Since $K$ is compact the map $\textnormal{id}_Y \times \phi : Y \times K \to Y \times K^\prime$ is closed (see \cite[3.6.3]{brw}). Let $q^\prime = (\textnormal{id}_Y \times \phi)\restriction X$ and $X^\prime = q^\prime(X)$. Since $X$ is closed in $Y \times K$, the set $X^\prime$ is closed in $Y \times K^\prime$. Besides, the map $q^\prime : X \to X^\prime$ is closed and hence an $\R$-quotient map. On the other hand, the map $r \times \textnormal{id}_{K^\prime} : Y \times K^\prime \to Y^\prime \times K^\prime$ is a continuous retraction. So, the second function will be  $q^{\prime\prime} = (r \times \textnormal{id}_{K^\prime}) \restriction X^\prime$. Put $X^{\prime\prime} = q^{\prime\prime}(X^\prime)$. Note that $q = q^{\prime\prime} \circ q^{\prime}$. Since $q(D)$ is dense in $q(X)$, we must have that  $q^{\prime\prime}(q^{\prime}(D))$ is dense in $X^{\prime\prime}$. Observe that $p_Y(D) \subseteq r(Y)$ implies that $q^{\prime\prime}(q^{\prime}(D)) = q^{\prime}(D)$. Hence, $q^{\prime}(D)$ is dense in $X^{\prime\prime}$. As a consequence $q^{\prime\prime}(X^\prime) = X^{\prime\prime} \subseteq \textnormal{cl}(q^{\prime}(D)) \subseteq \textnormal{cl}(q^{\prime}(X)) = \textnormal{cl}(X^{\prime}) = X^\prime$. It follows that $q^{\prime\prime}$ is a continuous retraction and hence an $\R$-quotient map. Therefore, $q = q^{\prime\prime} \circ q^{\prime}$ is also $\R$-quotient.
\endproof

Later on we will use the following notation. For a fixed cardinal $\kappa$ the space
$L_\kappa^\omega$ will be simply denoted by   $L$. Also we may assume that $L_\kappa$ has $\kappa + 1$ as a underlying set and $\kappa$ is the unique nonisolated point of $L_\kappa$. For each $G \in [\omega]^{<\omega}$, $u \in \kappa^G$ and $N \in [\kappa]^{\leq\omega}$ consider the clopen sets
$U_{G,N} = \{y \in L : y(G) \subseteq L_\kappa \setminus N\}$ and $V_{G,u} = \{y \in L : y \restriction_G = u\}$ of $L$. Put
$$
\mathcal{U}_L = \{U_{G,N} : G \in [\omega]^{<\omega}, N \in [\kappa]^{\leq\omega}\} \textnormal{ and } \mathcal{V}_L = \{V_{G,u} : G \in [\omega]^{<\omega}, u \in \kappa^G\}.
$$
Given $G \in [\omega]^{<\omega}$ and $V \in \mathcal{V}_L$ let $H_{G,V} := V \cap \{y \in L : y(G) \subseteq \{\kappa\}\}$. Observe that the family of clopen sets $\{U_{G,N} \cap V : N \in [\kappa]^{\leq\omega}\}$ is closed under countable intersections. Since $H_{G,V} = \bigcap\{U_{G,N} \cap V : N \in [\kappa]^{\leq\omega}\}$ and $L$ is Lindel\"of, it is easy to verify that
 the family $\{U_{G,N} \cap V : N \in [\kappa]^{\leq\omega}\}$ is a basis for the closed set $H_{G,V}$.

\begin{theorem}\label{CLTK}
If $K$ is compact and $X$ is a closed subspace of $L \times K$, then $X$ has a full $q$-skeleton.
\end{theorem}

\proof
For each $A \subseteq [X]^{\leq\omega}$ choose $\mathcal{S}(A) = \kappa \cap \bigcup\{p_n(p_L(A)) : n \in \omega\} \in [L_\kappa]^{\leq\omega}$ and define $r_A : L \to L$ by the rule $r_A(y)(n) = y(n)$ if $y(n) \in \mathcal{S}(A)$ and $r_A(y)(n) = \kappa$ otherwise. Note that $r_A$ is a continuous retraction and $r_A \restriction_{p_L(A)}$ is one-to-one. Besides it is easy to verify that $\mathcal{S}$ is $\omega$-monotone. Let $\mathcal{V}(A) = \{V_{G,u} : G \in [\omega]^{<\omega} \textnormal{ and } u \in \mathcal{S}(A)^G\}$.  Observe that $\mathcal{V}$ is also $\omega$-monotone.

\smallskip

For each $x,x^\prime \in X$ with $p_K(x) \not = p_K(x^\prime)$ fix a map $f_{x,x^\prime} \in C(K)$ such that $f_{x,x^\prime}(p_K(x)) \not= f_{x,x^\prime}(p_K(x^\prime))$. For every $A \in [X]^{\leq\omega}$ set $\mathcal{C}(A) = \{f_{x,x^\prime} : x,x^\prime \in A, p_K(x) \not = p_K(x^\prime)\} \in [C(K)]^{\leq\omega}$. It is clear that $\mathcal{C}$ is $\omega$-monotone. Let $\psi_A : \R^{C(K)} \to \R^{C(K)}$ the map defined by $\psi_A(\chi)(f) = \chi(f)$ if $f \in \mathcal{C}(A)$ and $\phi_A(\chi)(f) = 0$ otherwise, for each $\chi \in \R^{C(K)}$. Let $\phi_A  = \psi_A \circ \Delta C(K): K \to \R^{C(K)}$. Note that $\phi_A \restriction_{p_K(A)}$ is one-to-one. Observe that if $F \subseteq A \in [X]^{\leq\omega}$ and $W = [f_1,\ldots,f_n;B_1,\ldots,B_n] \in \mathcal{W}(\mathcal{C}(F))$, then $\phi_A^{-1}(W) = \phi_{F}^{-1}(W)$.

\smallskip

For $A \in [X]^{\leq\omega}$ we shall consider the map $q_A = r_A \times \phi_A \restriction_{X}$. Next, we proceed to  define the partially ordered set $\Gamma$ and the family $\{D_s : s \in \Gamma\}$ of subsets of $X$.

\smallskip

Pick $F \in [X]^{<\omega}$. For each  $G \in [\omega]^{<\omega}$, $V \in \mathcal{V}(F)$ and $W \in \mathcal{W}(\mathcal{C}(F))$ we define $N_{F,G,V,W} \in [\kappa]^{\leq\omega}$ and $x_{F,G,V,W} \in X$ as follows: If there exists a set $N \in [\kappa]^{\leq\omega}$ such that $X \cap (U_{G,N} \cap V \times \phi_F^{-1}(W)) = \emptyset$, then put $N_{F,G,V,W} \in [\kappa]^{\leq\omega} = N$, and  $N_{F,U,V,G} = \emptyset$ otherwise. If $X \cap (H_{G,V} \times \textnormal{cl}(\phi_F^{-1}(W))) \not= \emptyset$, then  fix an arbitrary point $x_{F,G,V,W}$ in this intersection, and choose $x_{F,G,V,W}$ arbitrarily in the other case. Consider the sets
$$
N_F = \bigcup\{N_{F,G,V,W} : G \in [\omega]^{<\omega}, V \in \mathcal{V}(F) \ \text{and} \ W \in \mathcal{W}(\mathcal{C}(F))\} \in [\kappa]^{\leq\omega}
$$
and
$$
A_F = \{x_{F,G,V,W} : G \in [\omega]^{<\omega}, V \in \mathcal{V}(F) \ \text{and} \ W \in \mathcal{W}(\mathcal{C}(F))\} \in [X]^{\leq\omega}.
$$
Finally choose $E_F \in [X]^{\leq\omega}$ so that $N_F \cap \mathcal{S}(X) \subseteq \mathcal{S}(E_F)$ and $A_F \subseteq E_F$.
For each $A \in [X]^{\leq\omega}$ let $\mathcal{E}(A) = \bigcup\{E_F : F \in [A]^{< \omega}\}$ and let $\mathcal{D}(A)$ be the closure of $A$ under $\mathcal{E}$. It is evident $\mathcal{E}$ and $\mathcal{D}$ are $\omega$-monotone.

\medskip

\textbf{Claim.} Let $A \in [X]^{\leq\omega}$. If $\mathcal{E}(A) \subseteq A$, then $q_A(A)$ is dense in $q_A(X)$.

\smallskip

{\it Proof of the Claim.} Fix $A \in [X]^{\leq\omega}$ such that $\mathcal{E}(A) \subseteq A$.  Let $O$ be a non-empty open set in $q_A(X)$. Note that we can choose $U \in \mathcal{U}_L$, $V \in \mathcal{V}(A)$, $W \in \mathcal{W}(\mathcal{C}(A))$ and $(y,z) \in X$ such that $q_A(y,z) = (r_A(y),\phi_A(z)) \in (U \cap V) \times W$ and $\big((U \cap V) \times \textnormal{cl}(W)\big) \cap q_A(X) \subseteq O$. We can suppose that $U = U_{G,N_0}$ for some $G \in [\omega]^{<\omega}$ and $N_0 \in [\kappa]^{\leq\omega}$ with $\mathcal{S}(A) \subseteq N_0$. Note that $r_A(y) \in H_{G,V}$. Write $A = \bigcup\{F_n : n \in \mathbb{N}\}$ where each $F_n$ is a finite set and $F_n \subseteq F_{n+1}$ for each $n \in \mathbb{N}$. Since $\mathcal{V}$, $\mathcal{W}$ and $\mathcal{C}$ are $\omega$-monotone, we can find $m \in \mathbb{N}$ such that  $V \in \mathcal{V}(F_m)$ and $W \in \mathcal{W}(\mathcal{C}(F_m))$. Put $F = F_m$. As it was pointed out above $\phi_A^{-1}(W) = \phi_{F}^{-1}(W)$.

We assert that $X \cap (H_{G,V} \times \textnormal{cl}(\phi_F^{-1}(W))) \not= \emptyset$. Assume on the contrary that $X \cap (H_{G,V} \times \textnormal{cl}(\phi_F^{-1}(W))) = \emptyset$. Then $H_{G,V} \cap p_L(X \cap (L \times \textnormal{cl}(\phi_F^{-1}(W)))) = \emptyset$. As $K$ is compact, the set $p_L(X \cap (L \times \textnormal{cl}(\phi_F^{-1}(W))))$ is closed. Since the family $\{U_{G,N} \cap V : N \in [\kappa]^{\leq\omega}\}$ is a basis for the closed set $H_{G,V}$, we can find $N \in [\kappa]^{\leq\omega}$ such that $(U_{G,N} \cap V) \cap p_L(X \cap (L \times \textnormal{cl}(\phi_F^{-1}(W)))) = \emptyset$. It follows that $X \cap (U_{G,N} \cap V \times \phi_{F}^{-1}(W)) = \emptyset$. By construction $N_F$ has been chosen in such a form that $X \cap (U_{G,N_{F}} \cap V \times \phi_{F}^{-1}(W)) = \emptyset$ and $N_F \cap \mathcal{S}(X) \subseteq \mathcal{S}(E_F) \subseteq \mathcal{S}(\mathcal{E}(A))  \subseteq \mathcal{S}(A) \subseteq N_0$. Using this fact and the fact that $r_A(y) \in H_{G,V}$, we may conclude that $y \in U_{G,N_F} \cap V$. As a consequence, we obtain that  $(y,z) \in X \cap \big((U_{G,N_{F}} \cap V) \times \phi_{F}^{-1}(W)\big)$, which is a contradiction. Hence. we can consider the point $x_0  = (y_0,z_0) = x_{F,G,V,W} \in X \cap (H_{G,V} \times \textnormal{cl}(\phi_F^{-1}(W)))$. In one hand, we have that $x_0 \in E_F \subseteq \mathcal{E}(A) \subseteq A$ and on the other hand we have that $y_0 \in H_{G,V}$ and $z_0 \in \textnormal{cl}(\phi_F^{-1}(W)) = \textnormal{cl}(\phi_A^{-1}(W))$. These two conditions imply that $q_A(x_0)  = (r_A(y_0),\phi_A(z_0)) \in (U \cap V \times \textnormal{cl}(W)) \cap q_A(X) \subseteq O$. This proves the Claim.

\medskip

Since $\mathcal{D}$ is $\omega$-monotone, the set $\Gamma = \mathcal{D}([X]^{\leq\omega})$ is up-directed and $\omega$-complete. Finally,  we define $D_A = A$ for each $A \in \Gamma$. If $A \in \Gamma$, then it follows from Claim 1 and Lemma \ref{LRQM}  that $q_A$ is an $\R$-quotient map. We shall prove that the families $\{q_A : A \in \Gamma\}$ and $\{D_A : A \in \Gamma\}$ form a full $q$-skeleton. In fact, (i) is a consequence of the Claim because of $\mathcal{E}(A) = A$, for each $A \in \Gamma$. To check (ii) assume that $A,B \in \Gamma$ and  $A \subseteq B$. Define $p_{B,A} := r_A \times \psi_A \restriction_{q_B(X)}$ and note that
$$
p_{B,A} \circ q_B = (r_A \times \psi_A) \circ (\hat{r}_B \times \phi_B) = (r_A \circ \hat{r}_B) \times (\psi_A \circ \phi_B) = r_A \times \phi_A = q_A.
 $$
The condition  $(iii)$ is immediate. Thus,   we conclude that the families $\{q_A : A \in \Gamma\}$ and $\{D_A : A \in \Gamma\}$ form a  $q$-skeleton. This $q$-skeleton is full because the conditions of Lemma \ref{LLRQS} are clearly satisfied.
\endproof

The following result was proved by I. Bandlow in \cite{ban} by using elementary submodels. Here, we obtain his result in a topological context.

\begin{corollary}
Let $K$ and $Z$ be compact spaces; suppose that $C_p(Z)$ is a continuous image of a closed subspace $X$ of $L \times K$. Then $Z$ is Corson compact.
\end{corollary}

\proof
Because of Theorem \ref{CLTK} the space $X$ has a full $q$-skeleton. Finally, it follows from Corollary \ref{CIC} that $Z$ is Corson.
\endproof


\end{document}